\documentclass[authoryear,preprint,review,2.5p,12pt]{elsarticle} 
\usepackage[utf8]{inputenc}
\usepackage{import}
\usepackage{siunitx}
\usepackage{amsmath,amsfonts,amssymb}
\usepackage{bookmark}
\usepackage{lscape}
\usepackage{graphicx}
\usepackage{setspace}
\usepackage{caption}
\usepackage{subcaption}
\usepackage{tabto}
\usepackage{algpseudocode}
\usepackage{algorithm}
\usepackage{booktabs}
\usepackage{soul}
\usepackage{psfrag}
\usepackage[T1]{fontenc}
\DeclareMathOperator*{\argmax}{arg\,max}

\usepackage{hyperref}
\makeatletter
\providecommand\HyPL@Entry[1]{}
\AddToHook{env/document/begin}{%
\@ifpackageloaded{preview}{
\ifPreview
 \let\Hy@FirstPageHook\relax
 \let\Hy@EveryPageAnchor\relax
\fi}{}}
\makeatother
\usepackage[version=4]{mhchem}
\setcitestyle{square,numbers,sort&compress,comma}
\usepackage{tikz}

\usepackage{threeparttable}
\usepackage{amsfonts} 
\usepackage{tikz}
\usepackage[normalem]{ulem}

\usepackage{amsmath}

\newcommand{\redline}{\raisebox{2pt}{\tikz{\draw[red,solid,line width = 1.2pt](0,0) -- (3mm,0)}}}

\newcommand{\grayline}{\raisebox{2pt}{\tikz{\draw[gray,solid,line width = 1.2pt](0,0) -- (3mm,0)}}}

\newcommand{\blackline}{\raisebox{2pt}{\tikz{\draw[black,solid,line width = 1.2pt](0,0) -- (3mm,0)}}}

\newcommand{\bluedashed}{\raisebox{2pt}{\tikz{\draw[blue,dashed,line width = 1.5pt](0,0) -- (3mm,0)}}}

\newcommand{\redcircle}{\raisebox{0.25pt}{\tikz{\filldraw [red] (0,0) circle (2pt)}}}

\newcommand{\bluecircle}{\raisebox{0.25pt}{\tikz{\filldraw [blue] (0,0) circle (2pt)}}}

\newcommand{\redsquare}{\raisebox{1pt}{\tikz{\filldraw [red] (0.1,0.1) rectangle (0.4,0.2)}}}

\newcommand{\bluesquare}{\raisebox{1pt}{\tikz{\filldraw [blue] (0.1,0.1) rectangle (0.4,0.2)}}}

\newcommand{\magsquare}{\raisebox{1pt}{\tikz{\filldraw [blue!50!red!50] (0.1,0.1) rectangle (0.4,0.2)}}}

\newcommand{\tr}[1]{\textcolor{black}{#1}}
\newcommand{\iter}[0]{N}

\title{BOATS: Bayesian Optimization for Active \tr{Control of} ThermoacousticS}
\date{November 2023}

\begin{document}
\begin{frontmatter}
\author{Bayu Dharmaputra\corref{cor1}}
\ead{bayud@ethz.ch}
\author{Pit Reckinger\corref{}}
\author{Bruno Schuermans\corref{}}

\author{Nicolas Noiray\corref{cor1}}
\ead{noirayn@ethz.ch}
\cortext[cor1]{Corresponding authors}
\address{CAPS Laboratory, Department of Mechanical and Process Engineering, ETH Z\unexpanded{\"u}rich, 8092, Z\unexpanded{\"u}rich, Switzerland}
\begin{abstract}

This investigation presents novel adaptive control algorithms specifically designed to address and mitigate thermoacoustic instabilities. Gas turbines are limited in their operational range due to thermoacoustic instability. Two control strategies are available to alleviate this issue: active and passive. Active control strategies have a wider flexibility than passive control strategies because they can adapt to the operating conditions of the gas turbine. However, optimizing the control parameters remains a challenge, especially if additional constraints have to be fulfilled, such as e.g. pollutant emission levels. To address this issue, we propose three adaptive control strategies based on Bayesian optimization. The first and foundational algorithm is the safeOpt algorithm, and the two adaptations that have been made are stageOpt and shrinkAlgo. The algorithms facilitate safe exploration within the control parameter space, ensuring compliance with the constraint function, while simultaneously optimizing the objective function. The Gaussian Process Regressor (GPR) is employed to approximate both the objective and constraint functions, with continuous updates occurring during iterations. The algorithms also enable the transfer of knowledge obtained from one operating point to another, thereby reducing the number of iterations needed to reach the optimal point. We demonstrate the effectiveness of the algorithms both numerically and through two distinct experimental validations. In the numerical demonstration, we employ a low-order thermoacoustic network model to simulate a single-stage combustor setup equipped with loudspeaker actuation and a gain-delay ($n-\tau$) controller for active stabilization. In the first experimental validation, we optimize the control parameters of a single-stage turbulent combustor with loudspeaker actuation and a gain-delay controller. For the second experimental validation, we apply the framework to a sequential combustor configuration utilizing nanosecond repetitively pulsed discharges (NRPD) as the control actuator. This demonstrates the framework's adaptability to various control actuation methods in turbulent combustors where control parameter optimization is required.

\end{abstract}

\end{frontmatter}

\section{Introduction}

Thermoacoustic instability is one of the major challenges in the development and application of modern gas turbines. It arises from a constructive interaction between acoustic perturbations and the heat release rate of the flame, which subsequently leads to large self-sustained acoustic pulsations. If it remains uncontrolled, it can lead to material failure and, therefore, expensive maintenance of gas turbines. 

Active and passive control strategies have been extensively studied in the context of thermoacoustic control. Due to their simplicity, passive control strategies have dominated the industrial gas turbine application. Passive control strategies can be grouped into two categories. The first category relies on the increase of the acoustic losses of the combustor \cite{Richards2003}. Typical examples in this category are Helmholtz dampers and quarter-wave resonators. The second category deals with reducing the acoustic driving of the combustion process.

Quarter-wave and Helmholtz resonators are the most commonly used devices, which belong to the first category.  Pandalai et al. \cite{Pandalai1998} showed an example of quarter-wave resonators placed in the cold section of the combustor upstream of the premixers in a GE aeroderivative engine. The authors noted that the resonator has logged more than 100000 hours of engine operations in factory testing and commercial operation. Quarter-wave resonators \tr{have the disadvantages of a typical narrowband response and a} length requirement \tr{which} can be prohibitive for typical frequencies in stationary engines \cite{Richards2000}. Bellucci et al. \cite{Bellucci2004} designed a Helmhotz damper that was later used in a Silo compressor (ALSTOM GT11N2) and also presented a nonlinear model to predict the loss coefficient and the natural frequency of the resonator. Significant research has been conducted \tr{in the recent years} to optimize the design of acoustic dampers for gas turbine combustors. For example,  Bothien et al\tr{.} \cite{bothien_2014} developed dampers made of interconnected cavities to widen their effective bandwidth and presented data from validation tests of their subwavelength damper concept in a heavy duty gas turbine. \tr{Bourquard and Noiray compared the volume and purge flow requirements of Helmholtz and quarter-wave resonators, and also showed that the optimum linear stability is achieved when these dampers are tuned to an exceptional point of the thermoacoustic system \cite{Bourquard2019288}. Later, } Miniero et al. \cite{Miniero2023} studied the effect \tr{modelled and investigated experimentally the generic problem} of \tr{periodic} hot gas ingestion into Helmholtz damper\tr{s} mounted on \tr{combustion chambers. They} \tr{presented an analytical model for robust  damper design, which can be used to limit the risk of passive control failure due to dynamic change of nonlinear damping and detuning}. 

Examples of the passive control strategy in the second category can be found in \cite{Lovett2002,Berenbrink2000}. \tr{In \cite{Lovett2002}, burners with} fuel injection \tr{at} different \tr{axial locations are used} to suppress acoustic pulsation\tr{. Such an axial staging leads to a bimodal distribution of the time delay corresponding to the convection of coherent equivalence ratio perturbations from the burner to the flame.} The dynamic phase converter presented by Noiray et al. in \cite{Noiray2009} is also based on an axial staging principle, but it relies on the bimodal time delay distribution of another type of convective coherent perturbations between the burner and the flames: the hydrodynamic ones. In fact, it works by converting \tr{long-wavelength} acoustic perturbations into \tr{short-wavelength} hydrodynamic perturbations \tr{thanks to} diaphragms \tr{in the injection channels}. The positions of the diaphragms are staggered to cancel out the fluctuating response of one-half of the \tr{reaction zone} with the other half. The stagger distance should be tuned to achieve effective cancelation at a target frequency. A more recent example is presented in \cite{AESOY2022} that is based on the alteration of flame transfer functions by modifying the geometry of the axial swirler and its position relative to the burner outlet. 

One of the disadvantages of passive control strategies is that they require extensive engine testing to acquire accurate knowledge of the system. For example, the instability frequency for different operating conditions must be known in order to design the appropriate damper geometry. Modifications to the burner geometry, either by adding diaphragms or changing the swirler geometry, might lead to additional unwanted pressure head loss. 

Active control strategies, on the other hand, could adapt to the changes in operating conditions. Seume et al. \cite{Seume1998} presented an active control strategy with pilot fuel flow modulation for the V84.3A gas turbine model. The control parameters consisted of gain and phase change ($n-\tau$ controller), and the actuators were multiple direct drive valves (DDV). The system was applied to a Siemens V94.3A heavy duty gas turbine and logged 18000 operating hours \cite{Hermann2005}. However, the bandwidth of the actuator is limited to 400~Hz. Another notable application of adaptive control in the GE aero derivative engine is presented in \cite{Pandalai1998}, where the combustor pulsation was measured with piezo sensors and the control actuation is achieved through splitting of the fuel between the inner, pilot, and outer rings of the burner. 

However, the implementation of active instability control (AIC) in commercial gas turbines is rather limited due to the lack of robust and cost-effective actuators \cite{Dharmaputra2023}. Nanosecond Repetitively Pulsed Discharges (NRPD) has shown to be an effective actuator to stabilize a sequential combustor as demonstrated in \cite{Dharmaputra2023}. Moeck et al. \cite{Moeck2013} implemented a feedback-based controller with NRPD to stabilize a swirling flame. One of the advantages of NRPD is that it does not require a significant modification of the combustor geometry. However, such an actuator would still require optimization of its control parameters to perform optimally and respecting some safety conditions. 

Active control algorithms have also been extensively studied for thermoacoustic stabilization. Gelbert et al. \cite{Gelbert2008} demonstrated for the first time the application of Model Predictive Control (MPC) to stabilize a turbulent swirl-stabilized combustor. The MPC requires a suitable model to work, in which the authors identified a priori. It was also noted that the system had to be in a stable state first, by adjusting their control parameters, for the MPC to work properly. There was no constraint implemented in both the states and the input variables; hence, the resulting Quadratic Programming (QP) problem can be solved in one step, enabling them to run the MPC algorithm with a sampling step of 1~ms. Adding input or state constraints might complicate optimization problems and increase the sampling time \cite{Schwenzer2021}. MPC algorithm can also be made robust with respect to some uncertainties; however, it always trades-off between robustness and optimality \cite{Schwenzer2021}. Nevertheless, the main drawback of model-based control is to find suitable models that can capture all the relevant operating conditions in gas turbines. The identification of the model parameters might be challenging and time consuming. Another control paradigm that addresses this issue is the adaptive control strategy. 

An example of adaptive control is the Extremum Seeking Controller (ESC). The demonstration of ESC for thermoacoustic stabilization was presented in \cite{Moeck2010}. In that study, the actuator is a loudspeaker equipped with a gain-delay controller. The ESC is used to adaptively tune the gain and delay parameters. ESC works by constantly perturbing the parameters with a low-amplitude sinusoidal function and subsequently computing the gradient of the objective function with respect to the control parameters. The gradient information is then used to drift the mean value of the parameters. The ESC algorithm was improved in \cite{Moeck2013} by encoding the slope information of the objective function so that the optimizer does not fall into a local maximum with zero gradient. However, there was no constraint encoded in the optimizer. 

A more recent adaptive control method for thermoacoustic stabilization with a loudspeaker is presented in \cite{Zhu2022}, where an active disturbance rejection control (ADRC) is implemented. The ADRC algorithm was originally introduced in \cite{Han2009} and works by treating the unknown plant dynamics as a disturbance that is tracked by an extended state observer (ESO). In principle, the ADRC algorithm is model-insensitive and can handle nonlinearities in the system. However, as demonstrated in \cite{Zhu2022}, tuning its free parameters would require proper modeling of the thermoacoustic system to obtain reasonable values. 

Most control synthesis methods, such as MPC, ARDC, and sliding mode controller (SMC) \cite{Tokat2015}, typically assume that all states are observable. If some states are not observable, an observer-based method can be employed such as Extender Kalman Filter (EKF). However, such an observer would require a model to work, and if the system is nonlinear, the convergence of the estimator might be an issue. Furthermore, the input action is typically assumed to be affine to the system, and can be written in the form of: 
\begin{equation}\label{eq:statespace}
    \dot{\mathbf{x}} = F(\mathbf{x}) + G(\mathbf{u}),
\end{equation}
where $\mathbf{x} \subset \mathbb{R}^n$ contains the state variables, $F: \mathbb{R}^n \rightarrow \mathbb{R}^n$ is some (non)-linear function of $x$ and $G: \mathbb{R}^m \rightarrow \mathbb{R}^n$ is some (non)-linear function of the input forcing $\mathbf{u} \subset \mathbb{R}^m$. 

Research in \cite{YU2023} has shown that the application of NRPD changes the flame response with respect to acoustic perturbations or, in other words, the flame transfer function (FTF). It is also demonstrated in \cite{Dharmaputra2023} that NRPD can stabilize a sequential combustor with continuous forcing and hence without a feedback loop. Therefore, these results hint that the NRPD forcing is not affine to the dynamics of the system. Therefore, most of the control synthesis methods are not directly implementable in this case. 

To address the issues, in this study, we propose an adaptive control method that is based on safe Bayesian optimization (safeOpt) which was first presented in \cite{Berkenkamp2016}. \tr{This} method is fundamentally data-driven, employing Gaussian Process Regressions (GPR) to approximate both objective and constraint functions. In that study, the safeOpt algorithm was used to optimize the parameters of a proportional derivative (PD) controller of a quadcopter. They showed that an optimum parameter combination could be found while also satisfying the safety constraint. The algorithm is further detailed in \cite{Berkenkamp2023}, where they showed that, by using context, the knowledge about good control parameters obtained at low tracking speeds can be transferred to fast tracking speeds. Khosavi et al. \cite{Khosravi2023} implement the same algorithm to optimize the gains of the PID cascade controller of a computer numerical control (CNC) grinding machine through both numerical and experimental tests. Their results showed that the algorithm performs 20$\%$ better than the nominal approach. Finally, it is worth mentioning the recent work of Reumschussel et al.\cite{Reumschussel2023}, who employed a Bayesian optimization-like strategy, albeit not based on the safeOpt algorithm, for experimental combustor design.

The present study introduces the safeOpt algorithm and two modified versions of the algorithm for active control of thermoacoustic instabilities in turbulent combustors. The algorithms are demonstrated in both numerical and experimental settings. The safeOpt algorithm will first be demonstrated in a numerical setup by employing a low-order thermoacoustic network model. Subsequently, the safeOpt algorithm and two additional modifications of it will be demonstrated in the experimental swirl-stabilized turbulent combustor setup. Finally, the safeOpt algorithm is demonstrated in a sequential combustor equipped with NRPD as the actuator.

To the best of our knowledge, the safeOpt algorithm has not been used in the domain of thermoacoustic control. The suitability of the safeOpt algorithm for thermoacoustic applications arises from its data-driven methodology (and hence does not require a model of the problem at hand), ability to adhere to constraints such as e.g. pollutant emission levels, input signal to actuators, turbine inlet temperature, etc., and straightforward implementation. Although the algorithm can incorporate information from a model into the prior mean function of the Gaussian Process Regression (GPR), our focus in this study is on scenarios where no such model is available, relying solely on measurements for optimization.

\section{Background Theory}\label{sec:background_theory}
A detailed explanation of the safeOpt algorithm is outlined in \cite{Berkenkamp2023}. Therefore, only a brief overview of the theoretical foundation is explained in this section. The problem statement is briefly outlined in Section \ref{sec:Prob_statement}. A brief overview of the Gaussian Process is discussed in Section \ref{sec:GP}. Afterward, the safeOpt algorithm and its modifications are explained in Section \ref{sec_safeOpt} and Section \ref{sec:mod_safeOpt}, respectively. Finally, the Bayesian context framework is discussed in Section \ref{sec:background_theory_bayesian_context}.

\subsection{Problem Statement} \label{sec:Prob_statement}
The goal of the proposed algorithm is to find control parameters $\mathbf{p}$ which optimize a scalar objective function $O$, under a certain constraint condition which is described by a constraint function $C$. The control parameters belong to a domain $\mathcal{P} \subset \mathbb{R}^{n}$. The objective function is defined as a map from the control parameter space to a scalar value: $O(\mathbf{p}) : \mathcal{P}\rightarrow\mathbb{R}$. Similarly, the constraint function is defined as: $C(\mathbf{p}) : \mathcal{P}\rightarrow\mathbb{R}$. Furthermore, we assume that there is an upper threshold value $T \in \mathbb{R}$ in which the system can be classified as safe: $C(\mathbf{p}) \leq T$. Note that the framework can be extended to include multiple constraint functions as described in \cite{Berkenkamp2023}. 

Both the objective function $O$ and the constraint function $C$ are not known a priori, but can be approximated by measurements for a given combination of control parameters $\mathbf{p}$. The algorithm will perform iteration updates and try to find the optimum point(s) of the aforementioned constrained optimization problem while ensuring that the safety condition is satisfied at each iteration $\iter$. Note that in our case we want to minimize the pressure pulsation; hence, this boils down to a minimization problem. Therefore, the optimization problem can be summarized as follows: 

\begin{equation}
    \min _{\mathbf{p} \in \mathcal{P}} O(\mathbf{p}) \text { subject to } C(\mathbf{p}) \leq T
\end{equation}

Since both the objective and constraint functions are not known a priori, an initial safe parameter set would need to be acquired. The set can be identified through simulations, expert domain knowledge, or some preliminary points evaluations. In our case, we perform multiple points evaluation in the domain $\mathcal{P}$ to obtain the initial set of safe parameters $\mathcal{S}_i \subset \tr{\mathcal{P}}$.

To expand knowledge of the safe parameter set beyond $\mathcal{S}_i$, the algorithm needs to infer whether some parameters $p^{*}$ that have not been evaluated are safe or unsafe. In this case, we make use of Gaussian Process (GP) model to approximate both the objective and constraint functions. Therefore, some regularity assumptions must be introduced for $O$ and $C$ \cite{Berkenkamp2023}. By using GP, we can construct a reliable confidence interval over $O$ and $C$, which allows us to satisfy the safety condition throughout the iterations with high probability. We denote $\hat{O}$ and $\hat
{C}$ as the GP approximation of $O$ and $C$, respectively. More precisely, since every measurement is contaminated by noise the GP approximations are defined as: 
\begin{equation} \label{eq:GPapproximation}
\begin{aligned}
        \hat{O}(\mathbf{p}) &=  O(\mathbf{p}) + \epsilon_o, \tr{\quad} \epsilon_o \sim \mathcal{N}(0,\sigma^2_o) \\
        \hat{C}(\mathbf{p}) &=  C(\mathbf{p}) + \epsilon_c, \tr{\quad} \epsilon_o \sim \mathcal{N}(0,\sigma^2_c) \\        
\end{aligned}      
\end{equation}
where $\epsilon_o$ and $\epsilon_c$ are gaussian random noise with zero mean and variance $\sigma_o$ and $\sigma_c$, respectively.

Note that due to the safety constraint condition, the algorithm might not be able to find the global optimum, however, it will aim to find the optimum parameters that are reachable from the initial safe set $\mathcal{S}_{i}$.

\subsection{Gaussian Process}\label{sec:GP}
In this work, GPs are used to approximate the objective funtion $O(\mathbf{p})$ and constraint function $C(\mathbf{p})$. GPs are non-parametric regression models which assume that the function values of the approximated function are random variables that have a joint Gaussian distribution \cite{Rasmussen}. 
A GP is described by a prior mean function and covariance function. The latter describes a covariance between two different parameter values: $\mathbf{p}, \mathbf{p'} \in \mathcal{P}$. A most commonly used term for the covariance function is kernel. In this study, we use a constant $K \in \mathbb{R}$ as the prior mean function and squared exponential kernels or Gaussian kernels as the covariance function. Note that the choice of the kernels is problem dependent, a detailed overview of possible kernels is available in \cite{Rasmussen}. The squared exponential kernel is defined as:

\begin{equation} 
     k(\mathbf{p},\mathbf{p'}) = \theta \exp\bigg(-\frac{d^2(\mathbf{p},\mathbf{p'})}{2}\bigg)
\end{equation}
\begin{equation} \label{eq:distancematrix}
     d^2(\mathbf{p},\mathbf{p'}) = (\mathbf{p}-\mathbf{p'})^{\intercal}\mathbf{L}^{-2}(\mathbf{p}-\mathbf{p'})
\end{equation}
where $\mathbf{L}$ is a diagonal matrix of positive real numbers representing the length scales: $\mathbf{L} = \mathrm{diag} (\mathbf{l})$, $\mathbf{l} \in \mathbb{R}^{n}_{+}$. Note that $n$ is the dimension of the control parameter space. The parameter $\theta$ represents the range of expected values of the difference in the value of the function and the previous mean function: $|K-O(\mathbf{p})|\leq 2\theta$ with 95\% probability. The length scales $\mathbf{l}$ represent how fast the covariance between neighbouring points decays with respect to their distance in the control parameters space $\mathcal{P}$. The last set of hyperparameters is the variance of measurement noise $\sigma^2_o$ and $\sigma^2_c$ in eq. \eqref{eq:GPapproximation}. In principle, the hyperparameters can also be optimized every time new data is acquired. However, as shown in \cite{bull2011}, this can lead to a poor result when using the maximum likelihood estimate to update the hyperparameters. Therefore, in the framework of safeOpt, the kernel hyperparameters are fixed from the beginning and treated as prior over functions. Hence, this represents the user's knowledge about the functions that are modeled. 

GPs can predict the function values $O(\mathbf{p^*})$ and $C(\mathbf{p^*})$, for any $\mathbf{p^*} \in \mathcal{P}$ based on the acquired data in the previous $\iter$ measurements. Until the end of this section, the notation of $\iter$ is shortened to $n$ for brevity. Conditioned on the measurements, the posterior distribution of the objective function (and equivalently the constraint function) is also Gaussian with the mean and variance as follows:
\begin{align}
    \mu_{n}(\mathbf{p}^*) &= \mathbf{k}_n(\mathbf{p}^*)(\mathbf{K}_n + \mathbf{I}_n \sigma_o^2)^{-1}\mathbf{\hat{O}}_n + K\\
    \sigma_n^2(\mathbf{p}^*) &= k(\mathbf{p}^*,\mathbf{p}^*)-\mathbf{k}_n(\mathbf{p}^*)(\mathbf{K}_n+\mathbf{I}_n\sigma_o^2)^{-1}\mathbf{k}_n^{\intercal}(a^*)
\end{align}    
where $\mathbf{\hat{O}}_n$ is the vector n observed values, $\sigma_o$ is the standard deviation of the observation noise, $\mathbf{k}_n(\mathbf{p}^*)$ is the covariance vector between the new point $\mathbf{p}^*$ and the observed data points, and $\mathbf{K}_n \in \mathbb{R}^{n\times n}$ is the covariance matrix of the observed data points, and $\mathbf{I}_n$ is an \tr{$n$ by $n$} identity matrix. It is worth mentioning that the kernel function $k$ for the objective and constraint functions may not have the same hyperparameters. For the remainder of the text we denote $k^o(\mathbf{p},\mathbf{p'})$, and $k^c(\mathbf{p},\mathbf{p'})$ as the kernel function for the objective and constraint function respectively.

\subsection{Safe Bayesian optimization (SafeOpt)}\label{sec_safeOpt}

\begin{figure}[t!]
\includegraphics[width=0.95\textwidth]{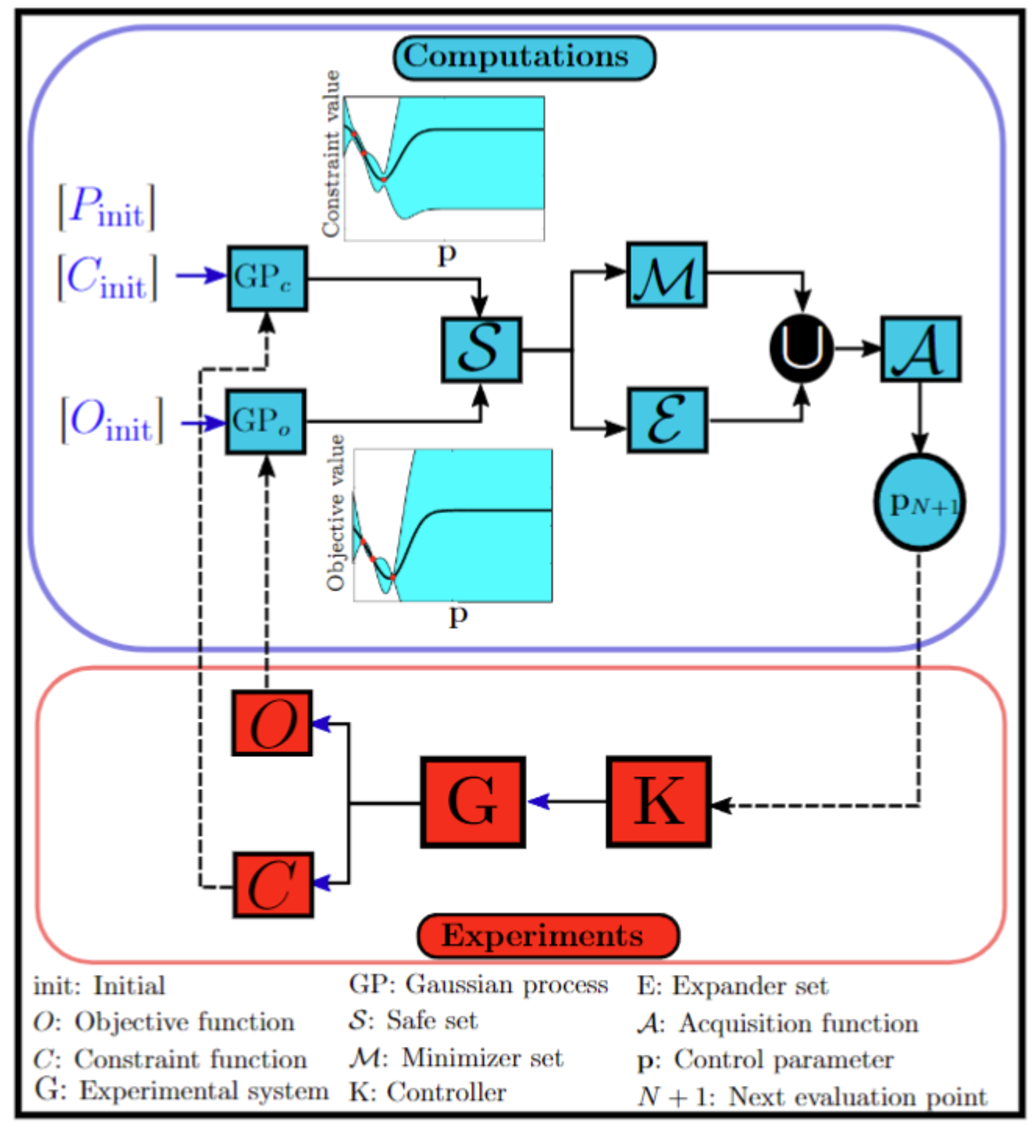}
    \caption[SafeOpt cartoon]{Visualization of the safeOpt algorithm: Gaussian process regressors model system responses at various control parameter values. Leveraging uncertainty bounds, the algorithm computes safe minimizer and expander sets. The union of these sets is input to an acquisition function, determining the next evaluation point. Objective and constraint function values are fed back to update the GPs.}
    \label{fig_safeOpt_vis}

\end{figure}

The SafeOpt algorithm, initially introduced in \cite{sui15} and later extended in \cite{Berkenkamp2023,Berkenkamp2016}, is visually represented in Figure \ref{fig_safeOpt_vis}. This algorithm operates by employing Gaussian Process Regressors (GPRs) to model the system's response in terms of objective and constraint function values. Utilizing uncertainty bounds, it calculates important sets, including the safe set $\mathcal{S} \subset \mathcal{P}$, the expander set $\mathcal{E} \subset \mathcal{P}$ and the minimizer set $\mathcal{M} \subset \mathcal{P}$. While fundamentally a Bayesian Optimization algorithm, safeOpt distinguishes itself by incorporating safety criteria throughout the entire iterative process. During each iteration, the algorithm seeks to identify the optimum point within the current safe set or expand the size of the safe set. This trade-off between exploration and exploitation is managed by selecting the point with the highest uncertainty in the objective function value. This ensures the adherence to additional safety criteria beyond optimizing the entire domain.

The safe set $\mathcal{S}$ is obtained by looking at the upper confidence bound of the GP estimate of the constraint function $U^c_n = \mu^c_n + 2\sigma^c_n$, and take the points which are below the threshold value $T$:

\begin{equation}\label{eq:safeset}
    \mathcal{S}_n = \mathcal{S}_{i}\cup \{\mathbf{p'}\in \mathcal{P}\| U^c_n(\mathbf{p'}) < T\}
\end{equation}
The choice of $2\sigma^c_n$ can be roughly interpreted as guaranteeing the safety with 95$\%$ probability per iteration. 

The potential minimizers which composes the minimizer set $\mathcal{E}$ is obtained by looking the upper bound of the objective function $U^o_n = \mu^o_n + 2 \sigma^o_n$ and the lower bound of the objective function $L^o_n = \mu^o_n - 2 \sigma^o_n$ which are inside the safe set $\mathcal{S}_n$ that satisfies the following: 

\begin{equation}\label{eq:minimizer}
    \mathcal{M}_n = \{\mathbf{p} \in \mathcal{S}_n | L^o_n(\mathbf{p}) < \min_{\mathbf{p}'\in \mathcal{S}_n} U^o_n(\mathbf{p}')\}
\end{equation}
which implies that the potential minimizers are the points in the safe set whose current lower bound estimate of the objective function is lower than the best upper bound.

Following \cite{Berkenkamp2016}, in order to define the expander set $\mathcal{E}_n$, an indicator function $e_n$ is first defined as follows:

\begin{equation}
    e_n(\mathbf{p}) = |\{\mathbf{p}' \in \overline{\mathcal{S}_n}| U^c_{n,(\mathbf{p},L^c_n(\mathbf{p}))}(\mathbf{p}') < T\}|,
\end{equation}
where $U^c_{n,(\mathbf{p},L^c_n(\mathbf{p}))}$ is the upper bound estimate of the constraint function based on the $n$ measurement points and an artificial measurement of $(\mathbf{p},L^c_n(\mathbf{p}))$. The indicator function $e_n({\mathbf{p}})$ computing the size of the previously unsafe sets $\overline{\mathcal{S}_n}$ that could potentially become safe if we hypothetically evaluated the point $\mathbf{p}$ and measured $L^c_n(\mathbf{p}))$ as the constraint function value. The expander set is therefore defined as:
\begin{equation}\label{eq:expander}
    \mathcal{E}_n = \{\mathbf{p} \in \mathcal{S}_n| e_n(\mathbf{p})>0\}.
\end{equation}
Essentially, the expander set comprises of the points that could potentially enlarge the current safe set. 

The next evaluation point $\mathbf{p}_{n+1}$ is acquired by taking the most uncertain point across the objective function which are inside the union of the potential minimizer set and the expander set:

\begin{equation}\label{eq:nextpoint}
    \mathbf{p}_{n+1} = \argmax_{\mathbf{p} \in \mathcal{E}_n \cup \mathcal{S}_n} U^o_n({\mathbf{p}}) - L^o_n({\mathbf{p}}). 
\end{equation}

The chosen acquisition function is widely known as "maximum uncertainty". The choice of the acquisition function will lead to a more exploratory behavior initially, which is caused by the fact that the most uncertain element typically lies on the boundaries of the safe region. Once the points close to the safety threshold $T$ are evaluated, the algorithm will evaluate the points in the potential minimizer and potential expanders alternatively. The proposed safeOpt algorithm is summarized in Algorithm \ref{alg:safeOpt}.

For demonstration purposes, we test the algorithm to perform an optimization in a one-dimensional space $\mathcal{P} \subset [0 ,10]$. The domain $\mathcal{P}$ is discretized with 200 points, an objective function with two minima and a constraint function with one minimum are chosen as follows: 

\begin{equation}
\begin{split}
    O &= 100(3\sin(2(p+1)^{0.8})-0.4p+7)\\
    C &= \frac{10}{(p+2)^{0.4}}+0.1(p-3)^{2}-4
\end{split} \label{eq:demo_1}
\end{equation}

\begin{figure}[b!]
\centering
\psfrag{-1000}[][]{\scriptsize -1000}
    \psfrag{0}[][]{\scriptsize 0~}
    \psfrag{0.5}[][]{\scriptsize 0.5~}
    \psfrag{-0.5}[][]{\scriptsize -0.5~}
    \psfrag{-1.5}[][]{\scriptsize ~~-1.5}
    \psfrag{1}[][]{\scriptsize 1~}
    \psfrag{1.5}[][]{\scriptsize 1.5~}
    \psfrag{2.5}[][]{\scriptsize 2.5~}   
    \psfrag{2}[][]{\scriptsize 2}
    \psfrag{4}[][]{\scriptsize 4}
    \psfrag{6}[][]{\scriptsize 6}
    \psfrag{8}[][]{\scriptsize 8}
    \psfrag{10}[][]{\scriptsize 10}  
    \psfrag{500}[][]{\scriptsize 500~}
    \psfrag{1000}[][]{\scriptsize 1000~}
    \psfrag{time (s)}[][]{\scriptsize time (s)}
    \psfrag{1500}[][]{\scriptsize 1500~}
    \psfrag{2000}[][]{\scriptsize 2000~}
    \psfrag{gain}[][]{\scriptsize $\mathbf{p}$}
    \psfrag{Iteration1}[][]{\scriptsize $\iter = 1$}
    \psfrag{Iteration5}[][]{\scriptsize $\iter = 5$}
    \psfrag{Iteration10}[][]{\scriptsize $\iter = 10$}
     \psfrag{Iteration20}[][]{\scriptsize $\iter = 20$}
    \psfrag{Iteration30}[][]{\scriptsize $\iter = 30$}
    \psfrag{Iteration40}[][]{\scriptsize $\iter = 40$}
    \psfrag{rmsPressure}[][t]{\scriptsize O (-)}
    \psfrag{rmsVoltage}[][t]{\scriptsize C (-)}
\includegraphics[trim=-0.5cm 0.05cm 0cm 0cm,clip,width=1\textwidth]{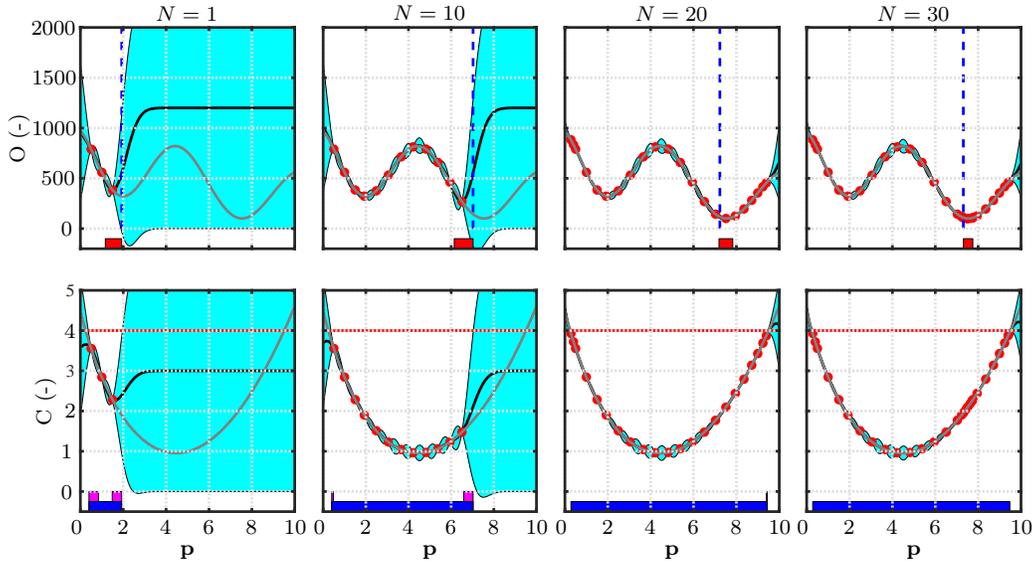}
    \caption[SafeOpt Introduction]{Demonstration Safe Optimization (safeOpt) algorithm after 
      1,10,20 and 30 iterations. (\bluesquare): safe set, (\redsquare): minimizer, (\magsquare): expander, (\blackline): mean prediction, (\bluedashed): next evaluation point, (\redcircle): measurement points, (\redline): safety constraint, (\grayline): real curve. The uncertainty of the prediction (2$\sigma$) is represented by the cyan shaded region.}
    \label{fig_safeOpt_demo}

\end{figure}

The performance of the algorithm is depicted in figure \ref{fig_safeOpt_demo}. Note that with the initial points for the algorithm are close to the first minimum point in the objective function which is not the global optimum. An optimization algorithm that relies on the local gradient of the function, such as the Extremum Seeking Controller (ESC), will easily be trapped in the first minimum. The algorithm spends the first 10 iterations to safely explore the parameter space. After about 20 iterations, the points with the constraint function value close to the threshold have been evaluated, and the algorithm starts evaluating a region around the global minimum within the domain.  

\begin{figure}[b!]
\centering
 \psfrag{-1000}[][]{\scriptsize -1000}
    \psfrag{0}[][]{\scriptsize 0~}
    \psfrag{0.5}[][]{\scriptsize 0.5~}
    \psfrag{-0.5}[][]{\scriptsize -0.5~}
    \psfrag{-1.5}[][]{\scriptsize ~~-1.5}
    \psfrag{1}[][]{\scriptsize 1~}
    \psfrag{1.5}[][]{\scriptsize 1.5~}
    \psfrag{2.5}[][]{\scriptsize 2.5~}   
    \psfrag{2}[][]{\scriptsize 2}
    \psfrag{4}[][]{\scriptsize 4}
    \psfrag{6}[][]{\scriptsize 6}
    \psfrag{8}[][]{\scriptsize 8}
     \psfrag{12}[][]{\scriptsize 12}
    \psfrag{16}[][]{\scriptsize 16}
    \psfrag{10}[][]{\scriptsize 10}  
    \psfrag{500}[][]{\scriptsize 500~}
    \psfrag{1000}[][]{\scriptsize 1000~}
    \psfrag{time (s)}[][]{\scriptsize time (s)}
    \psfrag{1500}[][]{\scriptsize 1500~}
    \psfrag{2000}[][]{\scriptsize 2000~}
    \psfrag{gain}[][]{\scriptsize $\mathbf{p}$}
    \psfrag{Iteration1}[][]{\scriptsize $\iter = 1$}
    \psfrag{Iteration5}[][]{\scriptsize $\iter = 5$}
    \psfrag{Iteration10}[][]{\scriptsize $\iter = 10$}
     \psfrag{Iteration20}[][]{\scriptsize $\iter = 20$}
    \psfrag{Iteration30}[][]{\scriptsize $\iter = 30$}
    \psfrag{Iteration40}[][]{\scriptsize $\iter = 40$}
    \psfrag{rmsPressure}[][t]{\scriptsize O (-)}
    \psfrag{rmsVoltage}[][t]{\scriptsize C (-)}
\includegraphics[trim=-0.3cm 0.05cm 0cm 0cm,clip,width=0.95\textwidth]{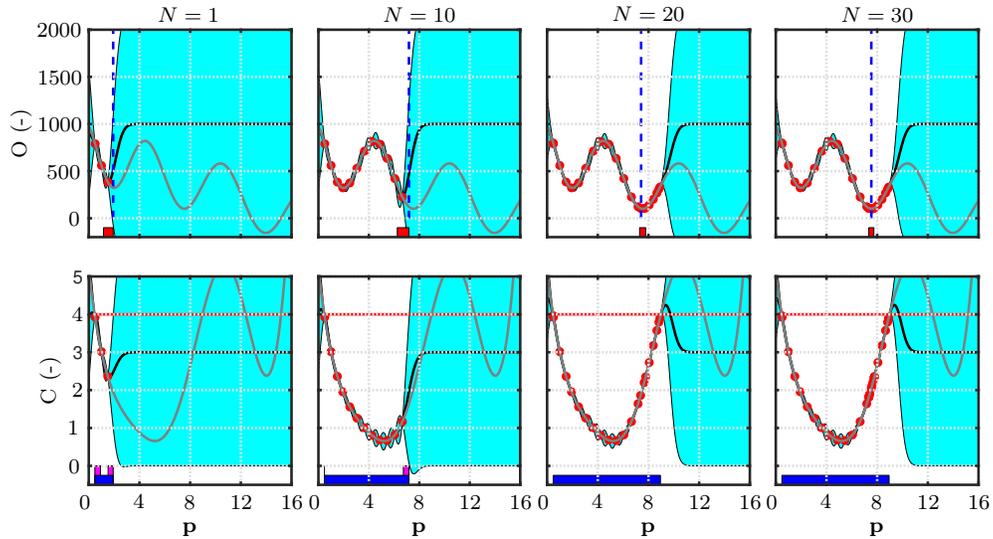}
   
    \caption[SafeOpt wells]{Demonstration Safe Optimization (safeOpt) algorithm with disjoint safe sets (see equation~\ref{eq:demo_2}) after 
      1,10,20 and 30 iterations. (\bluesquare): safe set, (\redsquare): minimizer, (\magsquare): expander, (\blackline): mean prediction, (\bluedashed): next evaluation point, (\redcircle): measurement points, (\redline): safety constraint, (\grayline): real curve. The uncertainty of the prediction (2$\sigma$) is depicted by the cyan shaded region.}
    \label{fig_safeOpt_demo_wells}

\end{figure}

Note that, as mentioned before, the algorithm will try to find the optimum location which is reachable from the initial safe set $\mathcal{S}_i$. Hence, the algorithm might encounter a problem when the safe set in the parameter space comprises of some disjoint sets and the initial safe set $\mathcal{S}_i$ does not contain some of the disjoint sets. To illustrate this issue, we keep the same objective function and modify the constraint function in equation~\ref{eq:demo_1} as follows: 

\begin{equation}
\begin{split}
    O &= 100(3\sin(2(p+1)^{0.8})-0.4p+7)\\
    C &= \frac{10}{(p+2)^{0.4}}+0.1(p-4)^{2}(1-0.7\sin(0.55p))-4
\end{split} \label{eq:demo_2}
\end{equation}
The parameter space is now from 0 to 16, $\mathcal{P} \subset [0, 16]$. The constraint function is slightly modified so that there are two disjoint safe sets in the whole domain. The performance of the safeOpt algorithm in this situation is shown in Figure~\ref{fig_safeOpt_demo_wells}. As can be seen, there are two safe regions in the domain $\mathcal{S} \subset {[0.5 ,9], [12.5 ,15]}$. Using the same initial safe set and hyperparameters, the algorithm is able to find the same minimum point safely. However, because the global minimum at $\mathbf{p} = 14$, is located in the second safe region, which is not reachable from the initial information in $\mathcal{S}_i$, the algorithm cannot reach it. However, if $\mathcal{S}_i$ included some points in the neighborhood of the global minimum, then the algorithm could have found the global minimum. Such a situation could indeed occur in an experimental setting; unless the algorithm is allowed to explore some unsafe point during the iterations or the initial safe set includes some points in all disjoint safe sets, the global minimum in the domain may not be uncovered. Nevertheless, the safeOpt algorithm will optimize the control parameter safely.

\begin{algorithm}[t!]
\caption{safeOpt Algorithm}\label{alg:safeOpt}
\begin{algorithmic}
\Require \\
~~-Control parameters domain $\mathcal{P}$\\
~~-GP kernel for the objective and constraint functions \\
~~-GP prior mean constants $K^o$ and $K^c$\\
~~-Safety threshold $T$\\
~~-Initial safe set $\mathcal{S}_i$ \Comment{In this study, it is generated by evaluating $N_{init}$ initial points}

\For{$\iter$ = 1,... $N_{max}$}
\State Find $\mathcal{S}_n$ (eq. \ref{eq:safeset})
\State Get potential minimizer set (eq. \ref{eq:minimizer})
\State Get possible expander set (eq. \ref{eq:expander})
\State Get the next evaluation point $\mathbf{p}_{n+1}$(eq. \ref{eq:nextpoint})
\State Obtain $\hat{O}(\mathbf{p}_{n+1})$ and $\hat{C}(\mathbf{p}_{n+1})$ through measurement
\State Update GPs with the new measurement data.
\EndFor \\
Select the best evaluated points $\mathbf{p}^*$

\end{algorithmic}
\end{algorithm}

\subsection{Modifications of safeOpt}\label{sec:mod_safeOpt}

We propose two types of modifications to adapt the algorithm to be more suited for thermoacoustic control. The modifications are named stageOpt and the shrinkAlgo. 
The stageOpt algorithm follows the one presented in \cite{SuiStageOpt2018}. It works by splitting the exploration and exploitation parts separately. In the current work, we let safeOpt algorithm work until $n_{s}$ iterations. SafeOpt superiority is used initially to minimize and explore the objective function safely; afterwards, the acquisition function in iteration $n_{s} + 1$ is switched to the minimum lower confidence bound (LCB):

\begin{equation}\label{eq:nextpointStageOpt}
    \mathbf{p}_{n+1} = \min_{\textbf{p} \in \mathcal{S}_i} \mu(\textbf{p})-2\sigma(\textbf{p})
\end{equation}

The second modification that we employ is to include a threshold in the objective function itself to give a second constraint. For the first $n_{s}$ iterations, safeOpt with a single constraint on the constraint function will be employed, afterwards a second constraint on the objective function $T_{o}$ is applied. The additional constraint would change the safe set and shrink its size. Therefore, we name it as "Shrinking" algorithm and abbreviate it as shrinkAlgo. More formally, the safe set after $n_{s}$ iterations is defined as: 

\begin{equation}\label{eq:safesetshrinking}
    \mathcal{S}_{n} = \{\mathbf{p'}\in \mathcal{P} | U^c_n(\mathbf{p'}) < T \cap U^o_n(\mathbf{p'}) < T_o\}.
\end{equation}
The motivation behind this is due to the fact that during the exploration phase, the algorithm could still evaluate points with high objective function values to enlarge the safe set. This situation could be undesirable if one wants to minimize the pressure pulsation as this will lead to the operating the system under high pulsation condition for long duration. Additionally, the computation of the expander set could be skipped to restrict the exploration of the points in the minimizer set, hence the choice of the next evaluation points can be written as follows: 

\begin{equation}\label{eq:nextpointshrinking}
    \mathbf{p}_{n+1} = \argmax_{\mathbf{p} \in \mathcal{S}_n} U^o_n({\mathbf{p}}) - L^o_n({\mathbf{p}}). 
\end{equation}
The choice of skipping the computation of the expander is optional, it could be done if the user is confidence that the current minimizer set could perform well and only fine exploration is required. The two proposed modifications are summarized in \tr{Algorithms} \ref{alg:stageOpt} and \ref{alg:shrinking}.

\begin{algorithm}[t!]
\caption{stageOpt Algorithm}\label{alg:stageOpt}
\begin{algorithmic}
\For{$\iter$ = 1,... $N_{s}$}
\State Follow algorithm \ref{alg:safeOpt}
\EndFor \\

\For{$\iter$ = $N_{s+1}$,... $N_{max}$}
\State Find $\mathcal{S}_n$ (eq. \ref{eq:safeset})
\State Get potential minimizer set (eq. \ref{eq:minimizer}) 
\State Get possible expander set (eq. \ref{eq:expander})
\State Get the next evaluation point $\mathbf{p}_{n+1}$(eq. \ref{eq:nextpointStageOpt})
\State Obtain $\hat{O}(\mathbf{p}_{n+1})$ and $\hat{C}(\mathbf{p}_{n+1})$ through measurement
\State Update GPs with the new measurement data.
\State 
\EndFor \\
Select the best evaluated points $\mathbf{p}^*$

\end{algorithmic}
\end{algorithm}

\begin{algorithm}[t!]
\caption{shrinkAlgo Algorithm}\label{alg:shrinking}
\begin{algorithmic}
\For{$\iter$ = 1,... $N_{s}$}
\State Follow algorithm \ref{alg:safeOpt}
\EndFor \\

\For{$\iter$ = $N_{s+1}$,... $N_{max}$}
\State Find $\mathcal{S}_n$ (eq. \ref{eq:safesetshrinking})\Comment{Safe set shrinks due to additional constraint}
\State Get potential minimizer set (eq. \ref{eq:minimizer})
\If{Use Expander is TRUE}
\State Get possible expander set (eq. \ref{eq:expander})
\State Get the next evaluation point $\mathbf{p}_{n+1}$(eq. \ref{eq:nextpoint})
\Else
\State Get the next evaluation point $\mathbf{p}_{n+1}$(eq. \ref{eq:nextpointshrinking})
\EndIf
\State Obtain $\hat{O}(\mathbf{p}_{n+1})$ and $\hat{C}(\mathbf{p}_{n+1})$ through measurement
\State Update GPs with the new measurement data.
\State 
\EndFor \\
Select the best evaluated points $\mathbf{p}^*$

\end{algorithmic}
\end{algorithm}

\subsection{Bayesian Context}\label{sec:background_theory_bayesian_context}

Bayesian context is a framework that allows us to model the dependency of the approximated functions with respect to additional external parameter(s) which are called context variables $\mathbf{z}$ \cite{KrauseOng2011}. The idea is to include the functional dependence and to keep it fixed when selecting the next points to evaluate \cite{Berkenkamp2023}. In the thermoacoustic context, this could be a small change in operating points such as fuel flow, air flow, and hydrogen blending level. Assuming that the frequency of the oscillation does not change significantly, the information from the previously optimized parameters in another condition can be transferred to the current one. This could speed up the optimization process as the previously information can be seen as measurement points with enlarged uncertainties. The dependence on the external parameter is modeled by creating a new kernel $k_c(\mathbf{z},\mathbf{z'})$ which will be multiplied by the kernel over the parameters defined in eq. (3). The total kernel is then defined as follows: 

\begin{equation}
    k([\mathbf{p},\mathbf{z})],[\mathbf{p'},\mathbf{z'})]) = k(\mathbf{p},\mathbf{p'})\times k_{\phi}(\mathbf{z},\mathbf{z'}).
\end{equation}
Hence, assuming one has evaluated the objective and constraint functions and approximated them with GP at the context variable $\mathbf{z}$, the information can be carried out to the next context variable $\mathbf{z'}$. The uncertainty will be enlarged depending on the kernel values of $k_{\phi}(\mathbf{z},\mathbf{z'})$.

\section{Experimental Setup}

\subsection{Single stage combustor}

 \begin{figure}[b!]
    \centering
    \includegraphics[trim=0.0cm 0.05cm -0.5cm 0cm,clip,width=1\textwidth]{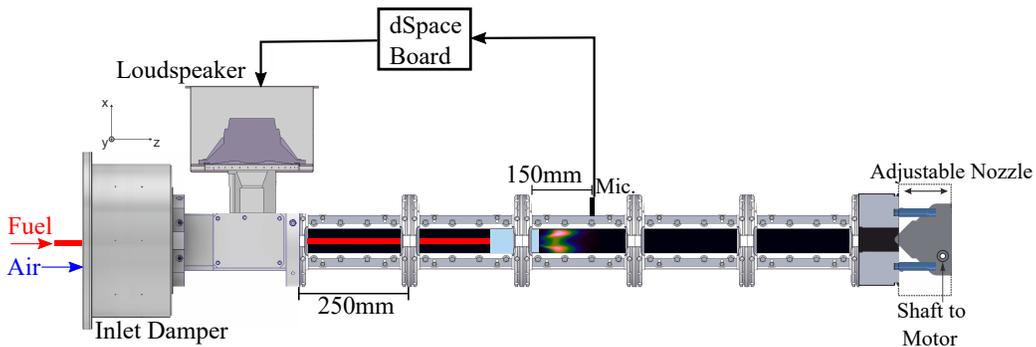}
    \psfrag{Loudspeaker}[][t]{\scriptsize Loudspeaker~~~~}
    \psfrag{Adjustable Nozzle}[][t]{\scriptsize Adjustable Nozzle}
    \psfrag{Fuel}[][t]{\scriptsize \textcolor{red}{Fuel}}
    \psfrag{Air}[][t]{\scriptsize \textcolor{blue}{Air}}
    \psfrag{250mm}[][]{\scriptsize 250 mm}
    \psfrag{Inlet Damper}[][]{\scriptsize Inlet Damper}
    \caption[ExpSetup]{Side view of the experimental setup. A loudspeaker is placed inside of an enclosure upstream of the burner.}
    \label{fig:ExpSetup}
\end{figure}

The cut view of the experimental setup can be seen in Figure \ref{fig:ExpSetup}. The experimental setup is similar to that in \cite{DHARMAPUTRAPROCIETF}. The experimental setup consists of an inlet plenum, adjustable inlet orifice, a loudspeaker, an axial swirler, and an adjustable piston at the end of the test rig. The ducts are made up of $250~\mathrm{mm}\times62~\mathrm{mm} \times62~\mathrm{mm}$ modules which are connected in series. The adjustable piston at the end of the test rig allows for a variation of the outlet orifice area. The presence of both inlet and outlet adjustable orifices allows the adjustment of the nominal thermoacoustic stability of the setup. A microphone is placed inside a water-cooled flush mounted plate in the combustion chamber module. In the technically premixed mode, the air is injected from the inlet plenum module whereas the fuel, which is a mixture of \tr{H$_2$} and \tr{CH$_4$}, is injected from the lance and delivered through 8 small holes downstream of the axial swirler inside the burner. 

The signal from the microphone is connected to a dSpace board (DS1104) where the controller is programmed to give an output voltage signal for the loudspeaker. In this study, the gain delay controller is employed; hence, the manipulation consists of delaying the microphone by $\tau$ milliseconds and multiplying the signal by a gain $n$. 

Three different operating conditions are considered, and they are enumerated as OP1, OP2, and OP3. The summary of important quantities for each operating condition is summarized in table \ref{table:OP_singlestage}. The instability frequencies for all operating conditions are around 200~Hz. 

\begin{table}[t!]
\begin{center}
    
\begin{tabular}{|l|l|l|l|l|l|}
\hline
OP & $\dot{m}_\mathrm{CH4}$ {[}g/s{]} & $\dot{m}_\mathrm{H_2}$ {[}g/s{]} & $\dot{m}_{air}$ {[}g/s{]} & $\phi$ {[}-{]} & $f_o$ {[}Hz{]}\\ \hline
1           & 0.57                                  & 0.02                                 & 15.35                              & 0.684                   & 200                     \\ \hline
2           & 0.64                                  & 0.01                                 & 15                                 & 0.753                   & 180                     \\ \hline
3           & 0.44                                  & 0.05                                 & 16                                 & 0.575                   & 230                     \\ \hline
\end{tabular}
\caption{Operating conditions for single stage combustor setup. $f_0$ denotes the instability frequency.} \label{table:OP_singlestage}
\end{center}
\end{table}

\subsection{Sequential combustor}

\begin{figure}[t!]
    \centering        
    \includegraphics[trim=0.0cm 0.05cm 0cm 0cm,clip,width=1\textwidth]{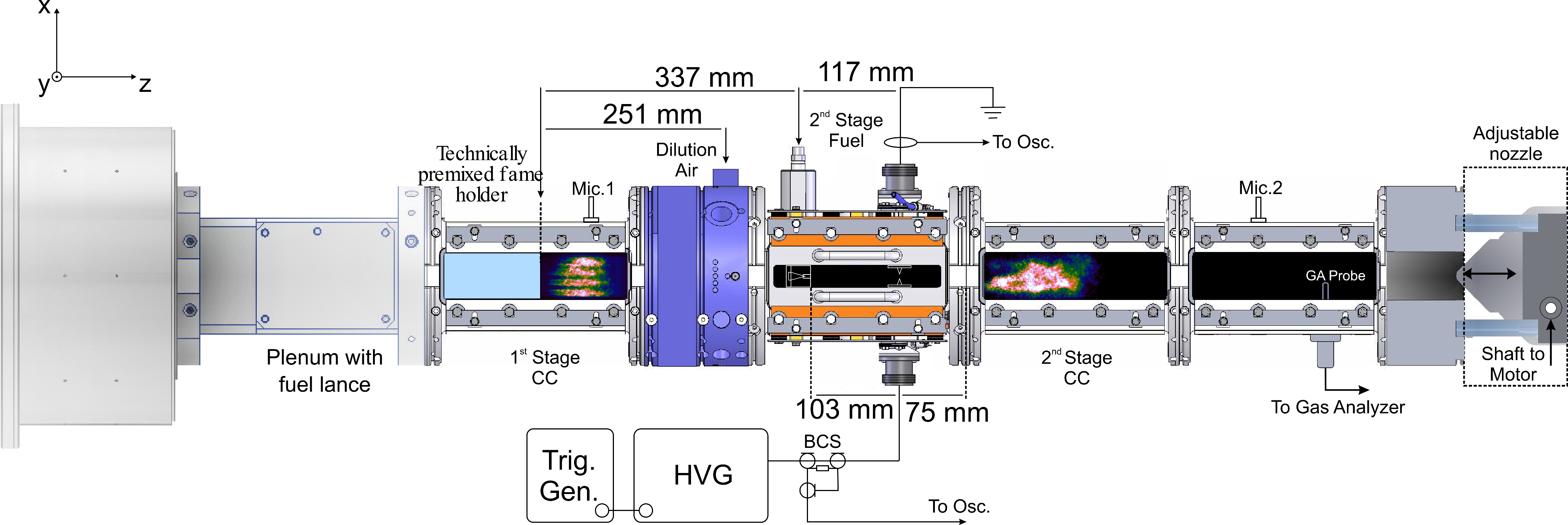}
    \caption[ExpSetupSeq]{Lab-scale sequential combustor test-rig. CC- Combustion chamber, HVG- High Voltage Generator, BCS- Back Current Shunt.}
    \label{fig:ExpSetupSequential}
\end{figure}

The lab-scale sequential combustor is depicted in figure \ref{fig:ExpSetupSequential}. The setup consists of a plenum, a 4 $\times$ 4 array of jet flames anchored on a matrix burner, a combustion chamber with a cross section of 62 $\times$ 62 $\mathrm{mm^2}$, a dilution air section, a sequential burner featuring a mixing channel with a cross section of 25 $\times$ 38 $\mathrm{mm^2}$, a sequential or second-stage combustion chamber equipped with a motor-driven adjustable outlet orifice. This variable outlet geometry enables an online tuning of the acoustic reflection coefficient, and thus an independent control of the thermoacoustic instabilities, which is key for validating the NRPD-based control. 

The first stage combustor is fed with a mixture of natural gas and air, with the air preheated to 230\textdegree C and supplied from the plenum. A piezo sensor is placed on a flush mounted plate to monitor the acoustic pressure inside the first stage and denoted as Mic.~1 in the figure. A massflow of 18~g/s of dilution air at 25~C is introduced from the dilution air port and mixes with the hot gases from the first stage. A mixture of  hydrogen and natural gas is injected into the sequential injector. The sequential injector features an X-shaped vortex generator to enhance the mixing process. The total thermal power of the two flames is 73.4~kW. A pin-to-pin electrode configuration, with an inter-electrode distance of 5~mm, is located 10.3~cm downstream from the sequential fuel injector, and a gas analyzer probe is placed at 45~cm from the outlet of the second-stage burner to monitor the NO emissions. Another piezo sensor is placed downstream of the sequential flame to monitor the acoustic pressure pulsation in the second combustion chamber. The exhaust gas analysis is conducted using an ABB EL3040 gas analyzer equipped with an Uras26 infrared photometer, operating at a sampling rate of 1~Hz. Automatic calibration with integrated calibration cells and sealing tests of the propagation line were performed before the start of each measuring set. The device exhibits a relative extended uncertainty of 7.9\% for NO measurements within the range of 0–200 mg/m$^3$.

\section{Thermoacoustic Network Model}

\begin{figure}[t!]
    \centering
    \begin{psfrags}
    
    \includegraphics[trim=0.0cm 0.0cm 0cm 0cm,clip,width=0.9\textwidth]{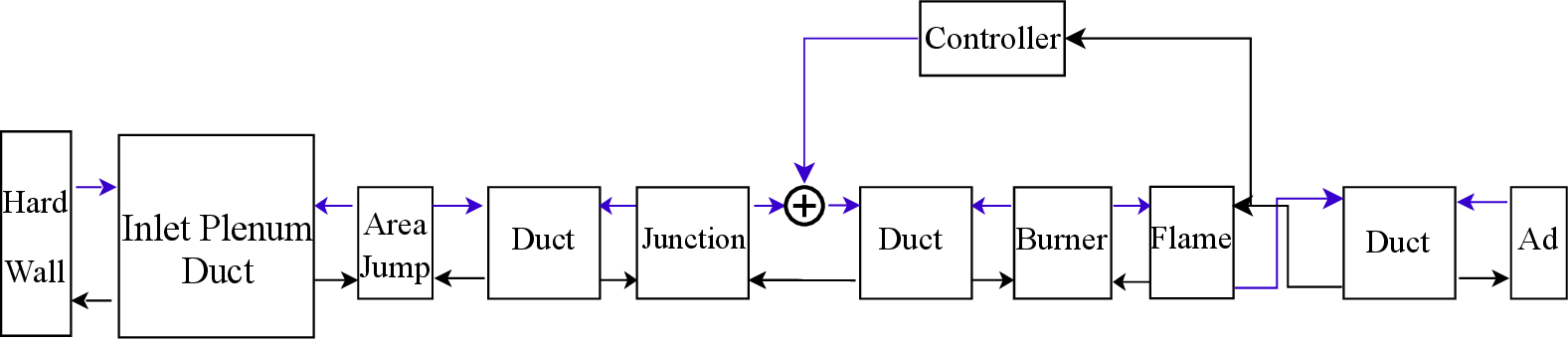}
    \caption{The thermoacoustic network model diagram. The blue arrow denotes the acoustic velocity perturbation $u'$, and the black arrow denotes the normalized acoustic pressure perturbation $p' /(\rho c)$. $\mathrm{Ad}$ is the downstream admittance boundary. The network model is implemented in MATLAB Simulink.}
    \label{fig:Lotandiagram}
    \end{psfrags}

\end{figure}

The safeOpt algorithm is first tested in a numerical setup, where the lab scale single stage combutor in figure \ref{fig:ExpSetup} is modeled with low order thermoacoustic network model. The network diagram is shown in figure \ref{fig:Lotandiagram}. The purpose of this study is not to model exactly each subsystem, the acoustic boundaries, flame transfer functions, and the pressure losses across the area expansions are tuned so that the system becomes unstable at around 200~Hz, which is the instability frequency of the setup. The network formulation follows the one presented in \cite{Bruno2010,BrunoThesis2003}. 

Each duct element is modeled as a perfect one-dimensional acoustic wave guide with velocity perturbations, $u'$, as input, and normalized acoustic pressure fluctuations $p'/(\rho c)$ as output. The acoustic perturbations can then be decomposed as forward- and backward-propagating waves, and expressed as follows:

\begin{equation}
    \frac{p'}{\rho c} = f + g
\end{equation}
\begin{equation}
    u' = f - g,
\end{equation}
where $\rho$ and $c$ are the gas density and the speed of sound, respectively.

The area jumps are modeled as compact area discontinuities with an equivalent length of $L_{eq}$ and a mean flow velocity in the orifice $\overline{U}_n$. By denoting $(.)_d$, $(.)_u$, $(.)_n$ as the acoustic quantities downstream, upstream of the area jump element, and inside the orifice, respectively, the governing equations can be written as: 

\begin{equation}
    A_d u'_d = A_u u'_u = A_n u'_n
\end{equation}
\begin{equation}
    \frac{p'_u - p'_d}{\rho c} - \frac{\overline{U_n}}{c}\zeta u'_n = \frac{L_{eq}}{c}\frac{du'_n}{dt}.
\end{equation}
The equations above are widely known as $L-\zeta$ model.

The flame is modeled as a compact element and isentropic assumption is employed for the components upstream and downstream of the flame. By using the Rankine-Hugoniot relation and do linearization to get the acoustic perturbations, the resulting flame transfer matrix can be written as:
\begin{equation}
    \begin{bmatrix}
 (\frac{p'}{\rho c})_d\\
 u'_d 
\end{bmatrix} =   \begin{bmatrix}
 \frac{(\rho c)_d}{(\rho c)_u} & 0 \\
 0 & 1+(\frac{T_d}{T_u}-1)FTF(\omega)
\end{bmatrix}    
\begin{bmatrix}
 (\frac{p'}{\rho c})_u\\
 u'_u 
\end{bmatrix}
\end{equation}
where $T_d$ and $T_u$ are the temperature downstream and upstream of the flame respectively. Note that the $T_{12}$ and $T_{21}$ elements above are equated to zero due to the low Mach number assumption. The Flame Transfer Function $FTF(\omega)$ is modeled with a first-order low-pass filter with a delay $\tau_f$: 
\begin{equation}
    FTF(\omega) = \frac{\exp(i\omega\tau_f)}{i\frac{\omega}{\omega_b} + 1},
\end{equation}
where $\omega = 2\pi f$ is the angular frequency, and $\omega_b$ the bandwidth of the low pass filter. For the time domain simulation, a tangent hyperbolic saturation function on the velocity perturbations upstream of the flame is utilized to saturate the amplitude of the downstream acoustic pressure fluctuations. 

The loudspeaker is modeled as a velocity perturbation located on the plenum side. Utilizing the pressure continuity in the junction element and the mass conservation, the relation between the acoustic quantities upstream and downstream of the junction is written as follows: 

\begin{equation}
    p'_u = p'_d    
\end{equation}
\begin{equation}
    u'_u = u'_d + u'_{LS},
\end{equation}
where $u'_{LS}$ is the acoustic velocity perturbations generated by the loudspeaker. Note that, the loudspeaker cavity and the electro-acoustic properties of the loudspeaker need to be taken into account if ones want to properly simulate the loudspeaker response as demonstrated in \cite{Lissek2011,RIVET2018}. However, it is not done in the numerical experiment performed in this section of the paper since we do not aim here at quantitatively reproducing the dynamics observed experimentally.

\begin{figure}[t!]
    \centering
\psfrag{50}[][]{\scriptsize 50}
    \psfrag{0}[][]{\scriptsize 0}
    \psfrag{-50}[][]{\scriptsize -50}
    \psfrag{100}[][]{\scriptsize 100}
    \psfrag{-100}[][]{\scriptsize -100}
    \psfrag{150}[][]{\scriptsize 150}
    \psfrag{-150}[][]{\scriptsize -150}
    \psfrag{200}[][]{\scriptsize 200}    
    \psfrag{300}[][]{\scriptsize 300}
    \psfrag{400}[][]{\scriptsize 400}
    \psfrag{500}[][]{\scriptsize 500}
    \psfrag{1}[][]{\scriptsize 1}
    \psfrag{2}[][]{\scriptsize 2}
    \psfrag{3}[][]{\scriptsize 3}
    \psfrag{4}[][]{\scriptsize 4}
    \psfrag{n}[][]{\scriptsize \tr{$n$}}
    \psfrag{FreqHz}[][]{\scriptsize $f$ (Hz)}
    \psfrag{GrowthRate}[][]{\scriptsize Growth Rate (1/s)}
    \psfrag{ngain}[][]{\scriptsize n}
    \includegraphics[trim=0cm 0.05cm 0cm 0cm,clip,width=0.5\textwidth]{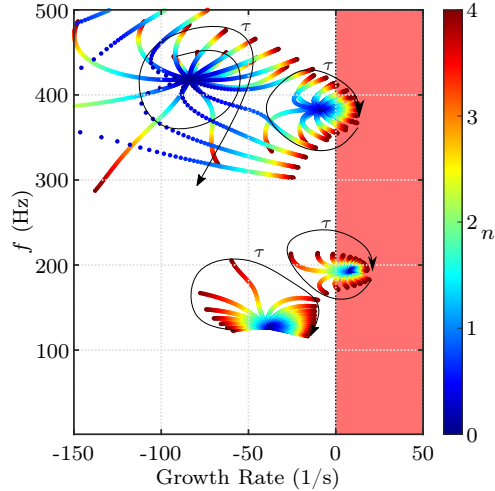}
    
    \caption{Maps of the eigenvalues of the system with varying gain ($n$) and the time delay ($\tau$) of the gain-delay controller. The red shaded region indicates the linearly unstable region.}
    \label{fig:eigval_maps}
\end{figure}

The controller takes the signal of the normalized acoustic pressure perturbations $p'/(\rho c)$ downstream of the flame and then delays it by $\tau$ milliseconds and multiplies it by a gain $n$, the signal is then placed into a saturation block, $L[a]=a$ for $|a|<5$ and $L[a]=\mathrm{sgn}(a)\times 5$ when $|a|\geq5$, so that the final output is bounded between $-5$ and 5. This saturation function represents the threshold voltage that we apply in the experiments to protect the loudspeakers from breaking due to high voltage values. The voltage output is then converted to $u'_{LS}$ by multiplying by a static gain, $K_{LS}$, of $-0.6$. More precisely, the voltage and velocity perturbations of the loudspeaker are expressed as: 

\begin{equation}
    V_{LS} = L\bigg[n \times\bigg(\frac{p'(t-\tau)}{\rho c}\bigg)_d \bigg]
\end{equation}
\begin{equation}
    u'_{LS} = K_{LS} V_{LS}
\end{equation}
It is worth mentioning that real loudspeakers will always have an effective frequency bandwidth in which the membrane will vibrate most efficiently. Since this behavior is not modeled in our case, all frequencies will pass through the loudspeaker without any attenuation. Hence it will be more probable in this numerical experiment that another acoustic mode is excited. 
By varying both parameters ($n$ and $\tau$), the system eigenvalues can be changed and will eventually shift all eigenvalues to the stable region. Figure \ref{fig:eigval_maps} shows the eigenvalues variation with respect to changes in both $n$ and $\tau$ of the controller. As seen, the unstable poles are around $f = 200~\mathrm{Hz}$ and 400~Hz. 

\section{Results}

In this section, the applications of the proposed algorithms are presented in both numerical and experimental settings. Section \ref{numerical} discusses the results obtained in the numerical setup, while Section \ref{experimental} shows the results from the experiments. In the numerical setup, only the safeOpt algorithm (Algorithm \ref{alg:safeOpt}) is employed, while, in the experimental setup with loudspeakers, all algorithms are employed. For both setups, optimization starts first with one parameter and then with two parameters. 

\subsection{Numerical validation}\label{numerical}

 \begin{figure}[t!]
    \centering
    \psfrag{-1000}[][]{\scriptsize -1000}
    \psfrag{0}[][]{\scriptsize 0~}
    \psfrag{0.5}[][]{\scriptsize 0.5~}
    \psfrag{-0.5}[][]{\scriptsize -0.5~}
    \psfrag{-1.5}[][]{\scriptsize ~~-1.5}
    \psfrag{1}[][]{\scriptsize 1~}
    \psfrag{1.5}[][]{\scriptsize 1.5~}
    \psfrag{2.5}[][]{\scriptsize 2.5~}   
    \psfrag{2}[][]{\scriptsize 2}
    \psfrag{4}[][]{\scriptsize 4~~}
    \psfrag{500}[][]{\scriptsize 500~}
    \psfrag{1000}[][]{\scriptsize 1000~}
    \psfrag{time (s)}[][]{\scriptsize time (s)}
    \psfrag{1500}[][]{\scriptsize 1500~}
    \psfrag{2000}[][]{\scriptsize 2000~}
    \psfrag{2500}[][]{\scriptsize 2500~}   
    \psfrag{gain}[][]{\scriptsize $n$~(-)}
    \psfrag{Iteration1}[][]{\scriptsize $\iter = 1$}
    \psfrag{Iteration5}[][]{\scriptsize $\iter = 5$}
    \psfrag{Iteration10}[][]{\scriptsize $\iter = 10$}
    \psfrag{Iteration20}[][]{\scriptsize $\iter = 20$}
    \psfrag{Iteration40}[][]{\scriptsize $\iter = 40$}
    \psfrag{rmsPressure}[][t]{\scriptsize O (Pa)}
    \psfrag{rmsVoltage}[][t]{\scriptsize C (V)}
    
    \includegraphics[trim=-0.5cm 0.05cm 0cm 0cm,clip,width=1\textwidth]{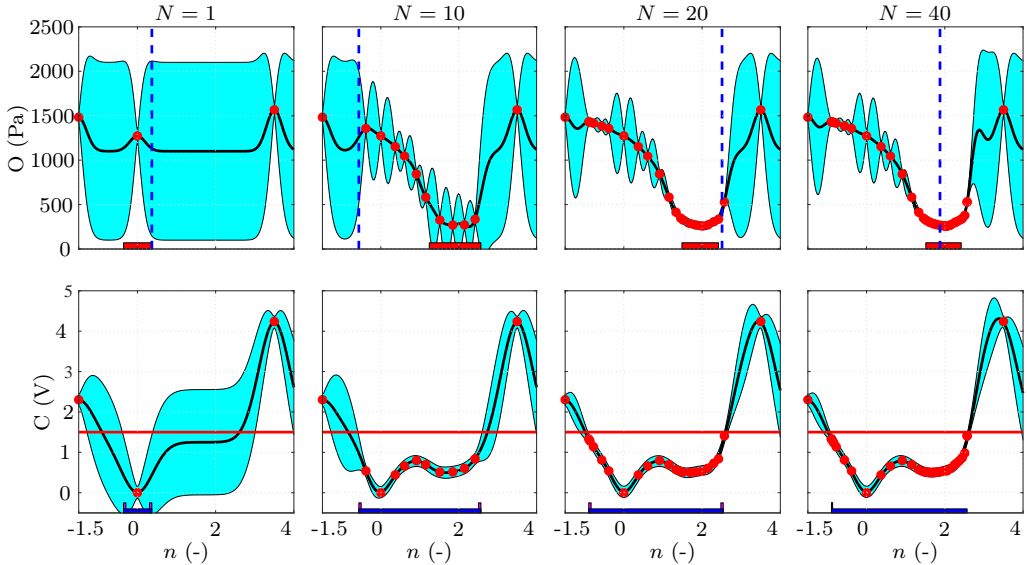}

    \caption[Gain Optimisation 1D safeopt simulation]{Controller $n$ optimization with SafeOpt after 1,10,20 and 40 iterations, the time delay, $\tau$, is fixed at 1.55~ms. (\bluesquare): safe set, (\redsquare): minimizer, (\magsquare): expander, (\blackline): mean prediction, (\bluedashed): next evaluation point, (\redcircle): measurement points, (\redline): safety constraint. The cyan-shaded region represents the uncertainty of the prediction (2$\sigma$). The global optimum of the objective function is already found within 10 iterations.}
    \label{fig:safeOpt_1D_Sim}

\end{figure}

\begin{figure}[t]
    \centering
        \psfrag{-1000}[][]{\scriptsize -1000}
        \psfrag{0}[][]{\scriptsize 0}
        \psfrag{1}[][]{\scriptsize 1}
        \psfrag{2}[][]{\scriptsize 2}
        \psfrag{3}[][]{\scriptsize 3}
        \psfrag{6}[][]{\scriptsize 6}
        \psfrag{-1.5}[][]{\scriptsize ~~-1.5}
        \psfrag{1.5}[][]{\scriptsize 1.5~}
        \psfrag{2.5}[][]{\scriptsize 2.5~}
        \psfrag{0.6}[][]{\scriptsize 0.6}
        \psfrag{1.2}[][]{\scriptsize 1.2}
        \psfrag{1.8}[][]{\scriptsize 1.8}
        \psfrag{80}[][]{\scriptsize 80}
        \psfrag{100}[][]{\scriptsize 100}
        \psfrag{120}[][]{\scriptsize 120}
        \psfrag{140}[][]{\scriptsize 140}
        \psfrag{160}[][]{\scriptsize 160}
        \psfrag{aaa}[][]{\scriptsize a}
        \psfrag{bbb}[][]{\scriptsize b}
        \psfrag{ccc}[][]{\scriptsize c}
        \psfrag{ddd}[][]{\scriptsize d}
        \psfrag{FreqHz}[][]{\scriptsize $f~(\mathrm{Hz})$}
        \psfrag{Uncontrolledxxx}[][]{\scriptsize No control~~~~}
        \psfrag{HighGainx}[][]{\scriptsize $n = 3.5$}
        \psfrag{BEP10xxx}[][]{\scriptsize $n = n^{\ast}_{10}$}
        \psfrag{BEP20xxx}[][]{\scriptsize $n = n^{\ast}_{20}$}
        \psfrag{200}[][]{\scriptsize 200}
        \psfrag{400}[][]{\scriptsize 400}
        \psfrag{600}[][]{\scriptsize 600}
        \psfrag{1000}[][]{\scriptsize 1000}
        \psfrag{3000}[][]{\scriptsize 3000}
        \psfrag{2000}[][]{\scriptsize 2000}
        \psfrag{-3000}[][]{\scriptsize -3000}
        \psfrag{-2000}[][]{\scriptsize -2000}
        \psfrag{0.5}[][]{\scriptsize 0.5}
        \psfrag{5}[][]{\scriptsize 5}
        \psfrag{times}[][]{\scriptsize $t$ (s)}
        \psfrag{gain}[][]{\scriptsize gain (-)}
        \psfrag{V}[][]{\scriptsize V}
        \psfrag{iter1}[][]{\scriptsize $\iter = 1$}
        \psfrag{iter17}[][]{\scriptsize $\iter = 17$}
        \psfrag{iter34}[][]{\scriptsize $\iter = 34$}
        \psfrag{PressurePa}[][t]{\scriptsize $p$ (Pa)}
        \psfrag{PSDdba}[][t]{\scriptsize $S_{pp}$ (dBa)}
        \psfrag{scaledpdf}[l][c][1][270]{\hspace{-2mm}\scriptsize $\hat{P}_p$}
        \psfrag{scaledpdf1}[l][c][1][270]{\hspace{-2mm}\scriptsize $\hat{P}_p$}
    \includegraphics[trim=-0.3cm 0.02cm 0cm 0cm,clip,width=1\textwidth]{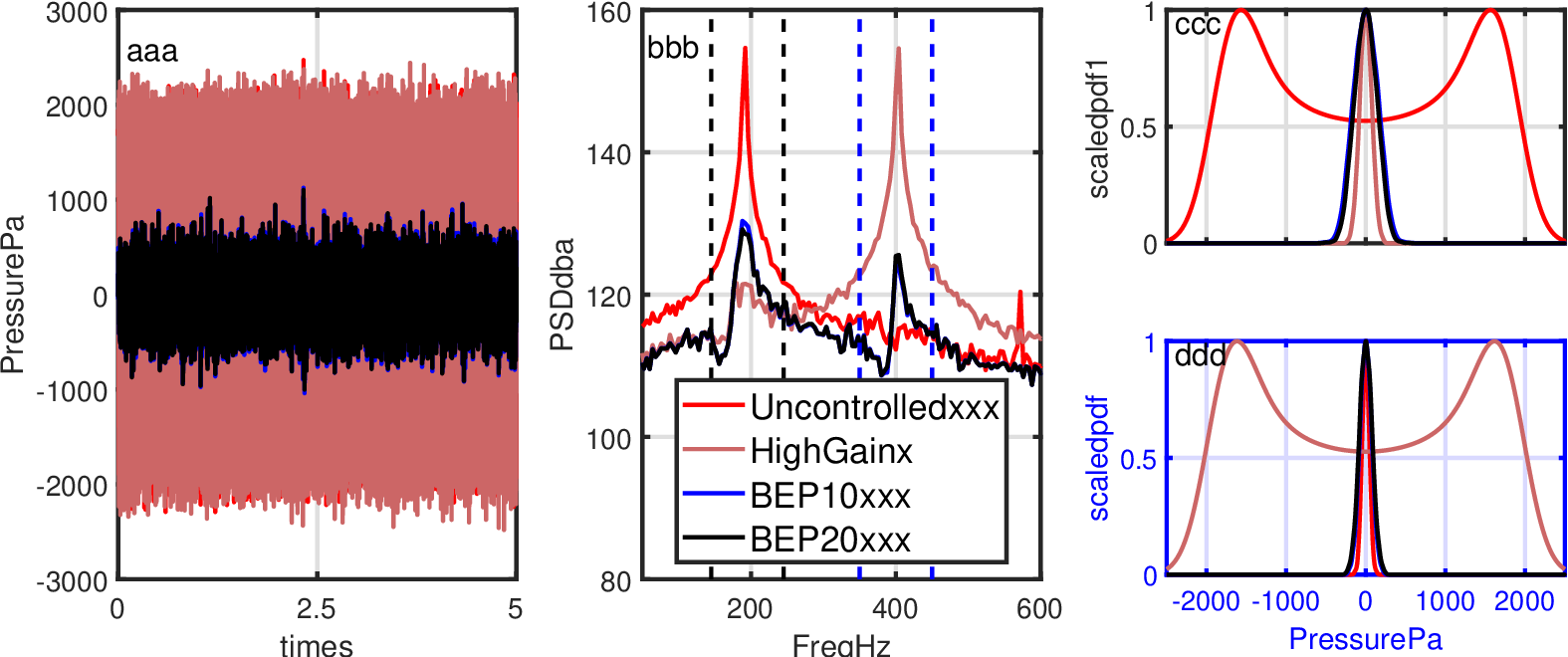}

    \caption[SafeOptimization 1D TT]{a) Simulated pressure time-trace and b) Power spectral density at four different gain values. $n = 3.5$ is the highest gain evaluated for the initial point. $n^{*}_{10}$ and $n^{*}_{20}$ denote the best evaluated $n$ after 10 and 20 iterations. c) and d) The pressure histogram of the bandpass filtered acoustic signal around 200 and 400 Hz respectively, the bandwidth of the filtering is indicated by the dashed line in b).} \label{fig:SafeOpt_sim_TT_HIST}

\end{figure}
For the gain-delay controller, the associated length scale $L$ in eq. \eqref{eq:distancematrix} can be casted into the following form: 

\begin{equation}
\mathbf{L} = \begin{bmatrix}
l_n & 0\\
0 & l_{\tau},
\end{bmatrix}
\end{equation}
where $l_n$ is the length scale of the gain parameter and $l_\tau$ is the length scale of the time delay parameter. 

\begin{table}[t!]
\begin{tabular}{@{}ccccc@{}}
\toprule
\multicolumn{1}{l}{\textbf{Hyperparameters}} & \textbf{\begin{tabular}[c]{@{}c@{}}$O$ \\ ($n$ Opt)\end{tabular}} & \textbf{\begin{tabular}[c]{@{}c@{}}$C$ \\ ($n$ Opt)\end{tabular}} & \textbf{\begin{tabular}[c]{@{}c@{}}$O$ \\ ($n-\tau$ Opt)\end{tabular}} & \textbf{\begin{tabular}[c]{@{}c@{}}$C$ \\ ($n-\tau$ Opt)\end{tabular}} \\ \midrule
$\theta$ {[}Pa - V{]}                        & 450                                                               & 0.65                                                              & 450                                                                    & 0.65                                                                   \\
$l_{n}$ {[}-{]}                              & 0.2                                                               & 0.4                                                               & 0.2                                                                    & 0.4                                                                    \\
$l_{\tau}$ {[}ms{]}                          & -                                             & -                                             & 0.75                                                                   & 0.75                                                                   \\
$\sigma$ {[}Pa - V{]}                        & 15                                                                & 0.05                                                              & 30                                                                     & 0.05                                                                  \\ \bottomrule
\end{tabular}
\caption{Hyperparameters of the Gaussian Process Regressors. Opt: optimization, $O$: objective function, $C$: constraint function}\label{tab:hyperparams}
\end{table}

First, only one parameter is optimized, which in this case is the gain $n$. The time delay of the controller is fixed at 1.55 ms. The objective function is the root mean square of the pressure pulsation after the flame, and the constraint function is the root mean square of the loudspeaker voltage. Both quantities are computed over a period of five seconds. Three initial points are fed to the algorithm, the initial points are $n = \{-1.5,0,3.5\}$, the initial safe set $\mathcal{S}_i$ is then obtained by fitting a Gaussian Process Regressor (GPR) to both the objective and the constraint function values. The domain of control parameters $\mathcal{P}$, which in this case contains only the $n$, is discretized by 100 uniform grid points from -1.5 to 4. The hyperparameters are listed in table \ref{tab:hyperparams}. The length scale and the prior variance of the kernel of the objective function $k^o(n,n')$ are set to 0.2 and 450 respectively. This can be interpreted that a distance of gain around 0.2-0.4 would yield completely different behavior, and the expected deviations from the mean value are 900 Pa. Whereas, the length scale and the prior variance of the kernel of the constraint function $k^c(n,n')$ are set to 0.4 and 0.65V, respectively. Noise variances $\sigma_o$ and $\sigma_c$ are set to 15 Pa and 0.05V, respectively. 

The gain optimization of the ($n$ - $\tau$) controller from the thermoacoustic network model is shown in Figure \ref{fig:safeOpt_1D_Sim}. As seen, the algorithm can safely find the global optimum in 10 iterations. After 20 iterations, more points close to the optimal location are evaluated. Additionally, the algorithm also tries to expand the safe set during the process. The size of the safe set, as shown by the blue rectangle in the bottom plot, increases throughout the iterations until the $20^{th}$ iteration. Subsequently, until the $40^{th}$ iterations, the algorithm evaluates almost exclusively the region close to the minimum. In the simulation, the modes at 200~Hz and at 400~Hz can be excited depending on the value of $n$. Figure \ref{fig:SafeOpt_sim_TT_HIST} shows the time trace, frequency spectra and the histogram of the filtered acoustic pressure pulsation around the instability frequencies. Without any control action, the system is unstable and the instability frequency is 200~Hz. When the gain is equal to 3.5, the higher mode at around 400~Hz becomes self excited. The best evaluated $n$ after 10 and 20 iterations, denoted by $n^{*}_{10}$ and $n^{*}_{20}$ exhibit almost the same performance in terms of pressure pulsation. This indicates that, in principle, 10 iterations are enough to optimize the gain safely. 

\begin{figure}[t]
    \centering
\psfrag{-1000}[][]{\scriptsize -1000}
        \psfrag{0}[][]{\scriptsize 0}
        \psfrag{1}[][]{\scriptsize 1}
        \psfrag{2}[][]{\scriptsize 2}
        \psfrag{3}[][]{\scriptsize 3}
        \psfrag{6}[][]{\scriptsize 6}
        \psfrag{-1.5}[][]{\scriptsize ~~-1.5}
        \psfrag{1.5}[][]{\scriptsize 1.5~}
        \psfrag{2.5}[][]{\scriptsize 2.5~}
        \psfrag{0.5}[][]{\scriptsize 0.5~}
        \psfrag{3.5}[][]{\scriptsize 3.5~}
        \psfrag{7}[][]{\scriptsize 7~}
        \psfrag{0.6}[][]{\scriptsize 0.6}
        \psfrag{1.2}[][]{\scriptsize 1.2}
        \psfrag{1.8}[][]{\scriptsize 1.8}
        \psfrag{1.25}[][]{\scriptsize 1.25}
        \psfrag{200}[][]{\scriptsize 200}
        \psfrag{600}[][]{\scriptsize 600}
        \psfrag{1000}[][]{\scriptsize 1000}
        \psfrag{Gain}[][]{\scriptsize $n~(-)$}
        \psfrag{1400}[][]{\scriptsize 1400}
        \psfrag{2000}[][]{\scriptsize 2000~}
        \psfrag{gain}[][]{\scriptsize gain (-)}
        \psfrag{Pa}[][]{\scriptsize Pa}
        \psfrag{V}[][]{\scriptsize V}
        \psfrag{iter1}[][]{\scriptsize $\iter = 1$}
        \psfrag{iter17}[][]{\scriptsize $\iter = 17$}
        \psfrag{iter34}[][]{\scriptsize $\iter = 34$}
        \psfrag{PhaseShift}[][t]{\scriptsize $\tau~(ms)$}
    \includegraphics[trim=-0.3cm 0.02cm 0cm 0cm,clip,width=0.94\textwidth]{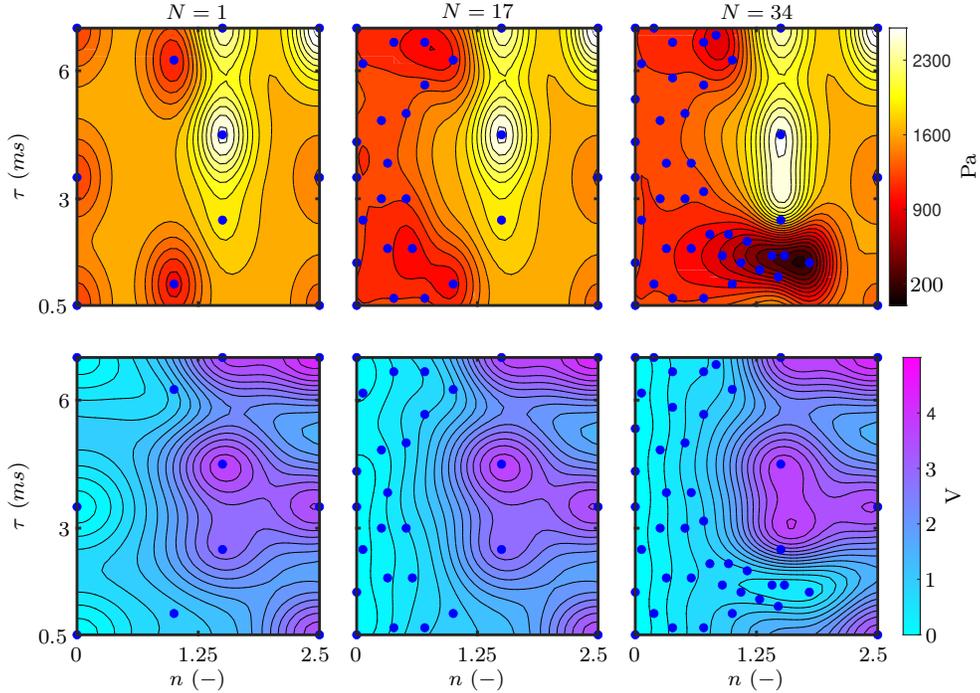}

    \caption[SafeOptimization 2D simulink]{$n-\tau$ optimization with safeOpt algorithm after 1, 17, and 34 iterations. The mean prediction of the surface map of the objective function (top) and constraint function (bottom). (\bluecircle): evaluation points.}
    \label{fig:SIM_n_tau_safeOpt}
\end{figure}

The optimization of the two control parameters optimization is shown in Figure \ref{fig:SIM_n_tau_safeOpt}. The domain of control parameters is a 50 $\times$ 50 uniform grid points with $n$ spanning from 0 to 2.5 and the $\tau$ ranges from 0.5~ms to 7~ms, $\mathcal{P} \subset \mathbb{R}^{[0~2.5]\times[0.5~7]}$. All kernel parameters are set the same as in the case with $n$ only optimization, with the addition of the length scale for $\tau$. The length scale for $\tau$ in $k^o(\mathbf{p},\mathbf{p'})$ is set to 0.4~ms, whereas for the $k^c(\mathbf{p},\mathbf{p'})$ they are set to 1~ms. Because the mode at around 400~Hz could be excited and the number of control parameters is now two, more initial points are required to initialize the algorithm for the computation of $\mathcal{S}_i$. Note that, as previously mentioned, the initial safe set could be given directly by the user if the user has knowledge of the control parameter space. \tr{Eleven} initial points are given to the algorithm, six points are on the left and right boundaries of the domain, and five points are around the middle of the domain. The number of iterations is set to 34 iterations which then amounts to 45 evaluated points in the domain. The algorithm spends the first 17 iterations to expand the safe set and essentially evaluate the left half-plane of the domain. After reaching 17 iterations, the algorithm starts to evaluate points with low pulsations. The best evaluated point after 34 iterations is $(n,\tau) = (1.8, 1.5)$. 

Despite the simplification of the thermoacoustic network model, the performance of the algorithm gives a preliminary indication that it could potentially work in the real system. The algorithm is applied to the experimental setup shown in figure \ref{fig:ExpSetup} in the next section.

\subsection{Experimental validation} \label{experimental}
\subsubsection{Single-stage combustor with loudspeaker actuation} \label{experimental_single_stage}

It is worth mentioning that, in contrast to the numerical simulation in Section \ref{numerical}, the experimental setup shown in figure \ref{fig:ExpSetup} only has one unstable acoustic mode over the whole control parameter space. Three operating conditions with different instability frequencies are considered and summarized in table \ref{table:OP_singlestage}. The three algorithms explained in Section \ref{sec:background_theory} are employed and compared with each other. Similarly to Section \ref{numerical}, the algorithms are first tested with one parameter optimization and continue with two parameters optimization. Additionally, the Bayesian context algorithm in Section \ref{sec:background_theory_bayesian_context} is applied to transfer the knowledge between different operating conditions to enhance the convergence speed of the algorithm. 

\begin{table}[t!]
\begin{tabular}{@{}ccccc@{}}
\toprule
\multicolumn{1}{l}{\textbf{Hyperparameters}} & \textbf{\begin{tabular}[c]{@{}c@{}}$O$ \\ ($n$ Opt)\end{tabular}} & \textbf{\begin{tabular}[c]{@{}c@{}}$C$ \\ ($n$ Opt)\end{tabular}} & \textbf{\begin{tabular}[c]{@{}c@{}}$O$ \\ ($n-\tau$ Opt)\end{tabular}} & \textbf{\begin{tabular}[c]{@{}c@{}}$C$ \\ ($n-\tau$ Opt)\end{tabular}} \\ \midrule
$\theta$ {[}Pa - V{]}                        & 450                                                               & 0.65                                                              & 450                                                                    & 0.65                                                                   \\
$l_{n}$ {[}-{]}                              & 0.2                                                               & 0.4                                                               & 0.2                                                                    & 0.4                                                                    \\
$l_{\tau}$ {[}ms{]}                          & -                                             & -                                             & 0.3                                                                   & 1                                                                   \\
$\sigma$ {[}Pa - V{]}                        & 15                                                                & 0.05                                                              & 30                                                                     & 0.075                                                                  \\ \bottomrule
\end{tabular}
\caption{Hyperparameters of the Gaussian Process Regressors. Opt: optimization, $O$: objective function, $C$: constraint function}\label{tab:hyperparamsEXP}
\end{table}

The hyperparameters for the Gaussian Process Regression (GPR) are detailed in Table \ref{tab:hyperparamsEXP}. Notably, all values closely align with those employed in the numerical test cases. A specific adjustment is made for $l_\tau$ in the constraint function, where it is now configured to be 1 ms. This adjustment is made because of the presence of a single unstable mode, allowing for an expectation of a larger correlation distance. Furthermore, $\sigma^c$ is set to 0.075 V intentionally to induce greater uncertainty in the measurements. This deliberate increase in uncertainty promotes a more conservative algorithmic behavior, thereby ensuring that the safety criterion is not violated.

\begin{figure}[t!]
    \centering
        \psfrag{-1000}[][]{\scriptsize -1000}
        \psfrag{0}[][]{\scriptsize 0~}
        \psfrag{0.5}[][]{\scriptsize 0.5~}
        \psfrag{-0.5}[][]{\scriptsize -0.5~}
        \psfrag{-1.5}[][]{\scriptsize ~~-1.5}
        \psfrag{1}[][]{\scriptsize 1~}
        \psfrag{1.5}[][]{\scriptsize 1.5~}
        \psfrag{2.5}[][]{\scriptsize 2.5~}
        \psfrag{2}[][]{\scriptsize 2}
        \psfrag{4}[][]{\scriptsize 4~~}
        \psfrag{500}[][]{\scriptsize 500~}
        \psfrag{1000}[][]{\scriptsize 1000~}
        \psfrag{time (s)}[][]{\scriptsize time (s)}
        \psfrag{1500}[][]{\scriptsize 1500~}
        \psfrag{2000}[][]{\scriptsize 2000~}
        \psfrag{gain}[][]{\scriptsize $n$~(-)}
        \psfrag{Iteration1}[][]{\scriptsize $\iter = 1$}
        \psfrag{Iteration5}[][]{\scriptsize $\iter = 5$}
        \psfrag{Iteration10}[][]{\scriptsize $\iter = 10$}
        \psfrag{Iteration20}[][]{\scriptsize $\iter = 20$}
        \psfrag{Iteration40}[][]{\scriptsize $\iter = 40$}
        \psfrag{rmsPressure}[][t]{\scriptsize O (Pa)}
        \psfrag{rmsVoltage}[][t]{\scriptsize C (V)}
    \includegraphics[trim=-0.3cm 0.05cm 0cm 0cm,clip,width=1\textwidth]{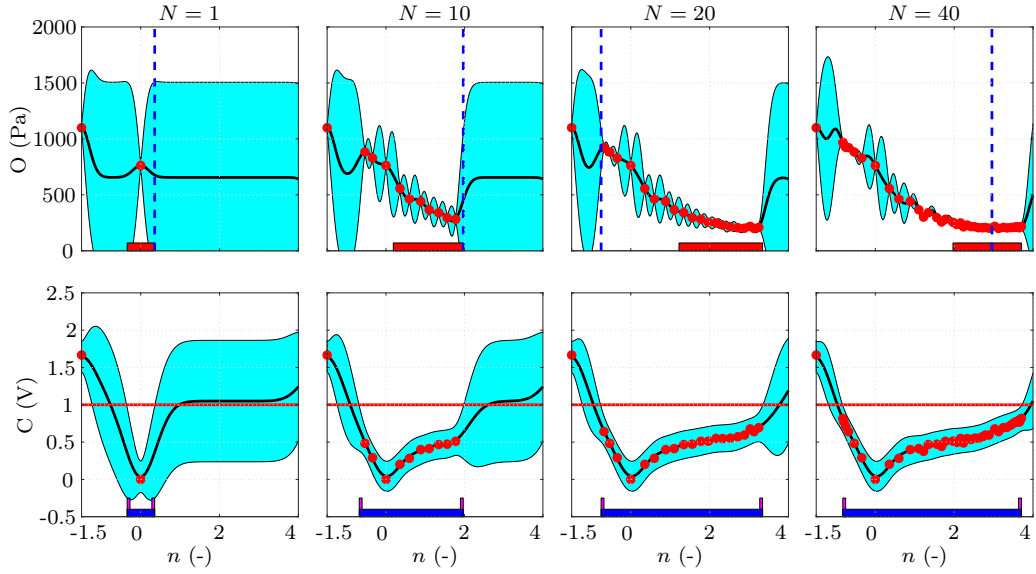}

    \caption[Gain Optimisation 1D safeopt]{Controller $n$ optimization at OP1 (see table \ref{table:OP_singlestage}) with SafeOpt algorithm after 1,10,20 and 40 iterations, the time delay, $\tau$, is fixed at 1.5~ms. (\bluesquare): safe set, (\redsquare): minimizer, (\magsquare): expander, (\blackline): mean prediction, (\bluedashed): next evaluation point, (\redcircle): measurement points, (\redline): safety constraint. The cyan-shaded region depicts the uncertainty of the prediction (2$\sigma$).}
    \label{fig:safeOpt_1D_Exp}
\end{figure}

The gain optimization of the gain-delay controller with safeOpt algorithm is shown in figure \ref{fig:safeOpt_1D_Exp}. Similarly to the numerical simulation, the control parameter space domain is discretized with 100 uniform grid points with n ranging from $-1.5$ to 4. Three initial points are given at $n = \{-1.5,0,4.5\}$. Note that in this case, one of the initial points is not inside the considered domain. As the purpose of the initial points is only to construct the initial safe set $\mathcal{S}_i$, this will not create any problem. The root mean square of the acoustic pulsation and the loudspeaker voltage are calculated by recording both signals for five seconds and then applying the rms operator. The safety threshold for the constraint function is $T~=~1\mathrm{V}$.  

\begin{figure}[t!]
    \centering
        \psfrag{-1000}[][]{\scriptsize -1000}
        \psfrag{0}[][]{\scriptsize 0}
        \psfrag{1}[][]{\scriptsize 1}
        \psfrag{2}[][]{\scriptsize 2}
        \psfrag{3}[][]{\scriptsize 3}
        \psfrag{6}[][]{\scriptsize 6}
        \psfrag{-1.5}[][]{\scriptsize ~~-1.5}
        \psfrag{1}[][]{\scriptsize 1~}
        \psfrag{1.5}[][]{\scriptsize 1.5~}
        \psfrag{2.5}[][]{\scriptsize 2.5~}
        \psfrag{0.6}[][]{\scriptsize 0.6}
        \psfrag{1.2}[][]{\scriptsize 1.2}
        \psfrag{1.8}[][]{\scriptsize 1.8}
        \psfrag{80}[][]{\scriptsize 80}
        \psfrag{100}[][]{\scriptsize 100}
        \psfrag{120}[][]{\scriptsize 120}
        \psfrag{140}[][]{\scriptsize 140}
        \psfrag{160}[][]{\scriptsize 160}
        \psfrag{aaa}[][]{\scriptsize a}
        \psfrag{bbb}[][]{\scriptsize b}
        \psfrag{ccc}[][]{\scriptsize c}
        \psfrag{ddd}[][]{\scriptsize d}
        \psfrag{FreqHz}[][]{\scriptsize $f~(\mathrm{Hz})$}
        \psfrag{Uncontrolledxx}[][]{\scriptsize No control~~~~}
        \psfrag{HighGainx}[][]{\scriptsize $n = 3.5$}
        \psfrag{BEP10}[][]{\scriptsize ~~$n = n^{\ast}_{10}$}
        \psfrag{BEP20}[][]{\scriptsize ~~$n = n^{\ast}_{20}$}
        \psfrag{200}[][]{\scriptsize 200}
        \psfrag{400}[][]{\scriptsize 400}
        \psfrag{600}[][]{\scriptsize 600}
        \psfrag{1000}[][]{\scriptsize 1000}
        \psfrag{3000}[][]{\scriptsize 3000}
        \psfrag{2000}[][]{\scriptsize 2000}
        \psfrag{-3000}[][]{\scriptsize -3000}
        \psfrag{-2000}[][]{\scriptsize -2000}
        \psfrag{0.5}[][]{\scriptsize 0.5}
        \psfrag{gain}[][]{\scriptsize gain (-)}
        \psfrag{V}[][]{\scriptsize V}
        \psfrag{95}[][]{\scriptsize 95}
        \psfrag{110}[][]{\scriptsize 110}
        \psfrag{125}[][]{\scriptsize 125}
        \psfrag{iter1}[][]{\scriptsize $\iter = 1$}
        \psfrag{iter17}[][]{\scriptsize $\iter = 17$}
        \psfrag{iter34}[][]{\scriptsize $\iter = 34$}
        \psfrag{PressurePa}[][t]{\scriptsize $p$ (Pa)}
        \psfrag{times}[][t]{\scriptsize $t$ (s)}
        \psfrag{PSDdba}[][t]{\scriptsize $S_{pp}$ (dBa)}
        \psfrag{scaledpdf}[l][c][1][270]{\hspace{-2mm}\scriptsize $\hat{P}_p$}
        \psfrag{scaledpdf1}[l][c][1][270]{\hspace{-2mm}\scriptsize $\hat{P}_p$}
    \includegraphics[trim=-0.3cm 0.02cm 0cm 0cm,clip,width=1\textwidth]{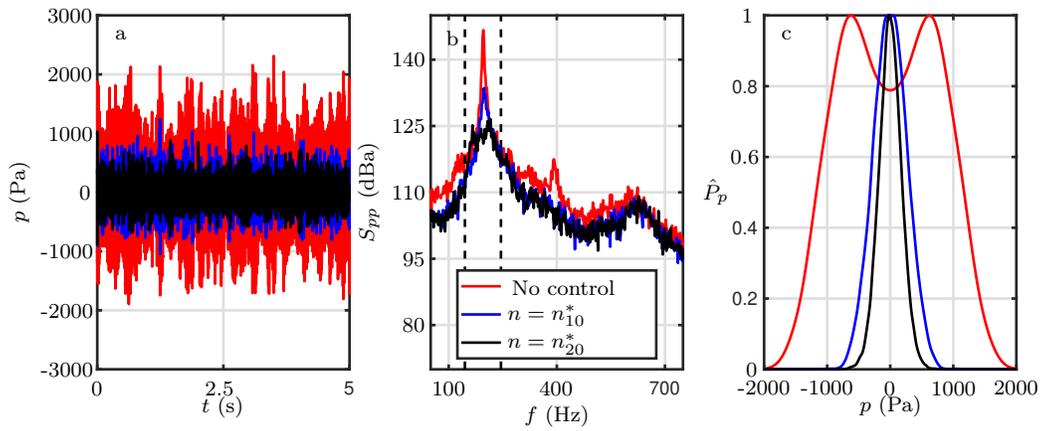}

    \caption[SafeOptimization 1D TT]{a) The pressure time-trace and b) Power spectral density at three different gain values. $n^{*}_{10}$ and $n^{*}_{20}$ denote the best evaluated $n$ after 10 and 20 iterations, respectively. c) The pressure histogram of the bandpass filtered acoustic signal around 200 Hz, the filtering bandwidth is indicated by the dashed line in b).}\label{fig:SafeOpt_exp_TT_HIST}
\end{figure}

Similar to the numerical simulation, the algorithm initially spends the first 10 iterations expanding the safe set, as evident from the growth in its size. After these initial 10 iterations, the algorithm shifts its focus to evaluating points with low pulsation. It is worth noting that, in contrast to the numerical simulation, when the gain values fall within the range of 1.5 to 3.5, the resulting pressure root mean square (rms) values are nearly identical. Figure \ref{fig:SafeOpt_exp_TT_HIST} provides visual representations of the time trace, frequency spectra, and the scaled probability density function of the acoustic pressure signal. Specifically, Figure \ref{fig:SafeOpt_exp_TT_HIST}c illustrates that the system stabilizes after just 10 iterations with the best-evaluated value of \tr{``$n$"}. However, after 20 iterations, a more favorable \tr{``$n$"} in terms of pressure rms is found.

As the iterations progress beyond 20, there is no improvement in pressure rms. The algorithm continues to evaluate different values of \tr{``$n$"}, but the pressure rms remains the same, as indicated in Figure \ref{fig:safeOpt_1D_Exp}.

The inherent nature of the safeOpt algorithm involves a continuous trade-off between exploring or expanding the safe set and minimizing the objective function throughout all iterations. Consequently, the algorithm occasionally evaluates points with high pressure rms values, even when the optimal point is unlikely to be found in this region. For example, in Figure \ref{fig:safeOpt_1D_Exp}, at the $20^{th}$ iteration, the algorithm chooses a safe value of n which is likely to have a high objective function value. This characteristic can be advantageous in escaping local optima if they exist. However, it may be undesirable when such local optima are absent, resulting in the combustor operating with high pulsation with no tangible benefits. Therefore, to solve this issue, stageOpt and the shrinking algorithm, which are explained in Algorithm \ref{alg:stageOpt}, and Algorithm \ref{alg:shrinking}, respectively, are used.

 \begin{figure}[t!]
    \centering
    \psfrag{-1000}[][]{\scriptsize -1000}
    \psfrag{0}[][]{\scriptsize 0~}
    \psfrag{0.5}[][]{\scriptsize 0.5~}
    \psfrag{-0.5}[][]{\scriptsize -0.5~}
    \psfrag{-1.5}[][]{\scriptsize ~~-1.5}
    \psfrag{1}[][]{\scriptsize 1~}
    \psfrag{1.5}[][]{\scriptsize 1.5~}
    \psfrag{2.5}[][]{\scriptsize 2.5~}
    \psfrag{2}[][]{\scriptsize 2}
    \psfrag{4}[][]{\scriptsize 4~~}
    \psfrag{500}[][]{\scriptsize 500~}
    \psfrag{1000}[][]{\scriptsize 1000~}
    \psfrag{time (s)}[][]{\scriptsize time (s)}
    \psfrag{1500}[][]{\scriptsize 1500~}
    \psfrag{2000}[][]{\scriptsize 2000~}
    \psfrag{gain}[][]{\scriptsize gain (-)}
    \psfrag{Iteration1}[][]{\scriptsize $\iter = 1$}
    \psfrag{Iteration5}[][]{\scriptsize $\iter = 5$}
    \psfrag{Iteration10}[][]{\scriptsize $\iter = 10$}
    \psfrag{Iteration40}[][]{\scriptsize $\iter = 40$}
    \psfrag{rmsPressure}[][t]{\scriptsize O (Pa)}
    \psfrag{rmsVoltage}[][t]{\scriptsize C (V)}
    \psfrag{gain}[][t]{\scriptsize $n$}
    
    \includegraphics[trim=0cm 0.05cm 0cm 0cm,clip,width=1\textwidth]{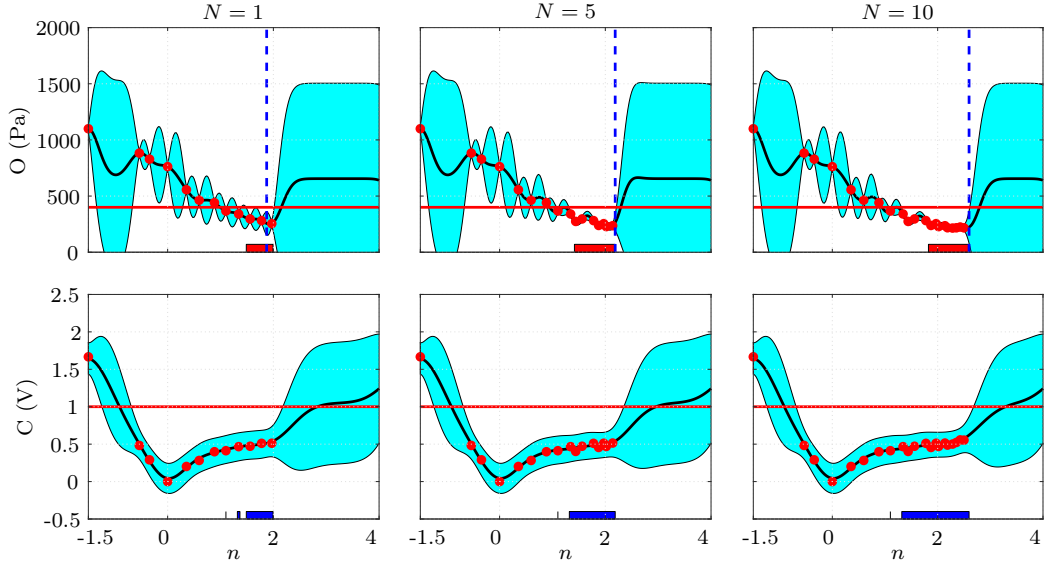}

    \caption[Gain Optimisation Shrinking]{Controller $n$ optimization at OP1 (see table \ref{table:OP_singlestage}) with shrinkAlgo after 1, 5, and 10 iterations. \tr{Note that $\iter$ displayed has been restarted to 1 once the additional constrain was activated after the 11th iteration.} The time delay, $\tau$, is fixed at 1.5~ms. (\bluesquare): safe set (the safe set fulfills both constraints on $O$ and $C$), (\redsquare): minimizer, (\blackline): mean prediction, (\bluedashed): next evaluation point, (\redcircle): measurement points, (\redline): safety constraint. Expander computation is excluded. The shrinking algorithm was activated after 10 iterations of SafeOpt. The safe set immediately shrinks after a secondary constraint on the objective function is applied.}
    \label{fig:n_shrinking}
\end{figure}

Figure \tr{\ref{fig:n_shrinking}} shows the results of the gain optimization with the shrinking algorithm. Following Algorithm \ref{alg:shrinking}, the first 10 iterations employ the regular safeOpt algorithm, at \tr{iteration} $N_s = 11$, an additional constraint $T_o$ is added to the objective function itself with a value of 450~Pa. Note that, in Figure \ref{fig:n_shrinking}, $\iter$ is restarted to 1 once the additional constraint is applied. As seen, the additional constraint leads to a shrinkage of the safe set. Additionally, the expander is not used; however, as iterations progress, the safety set is still growing because the acquisition function is the maximum uncertainty as described in eq \ref{eq:nextpointshrinking}. In this case, the benefit is clear, the algorithm does not evaluate points with high pressure pulsations since these points are now classified as unsafe. Note that, if the additional constraint is applied too early in the iterations, it may not be able to find a safety set, as the points with an rms pressure value below the threshold $T_o$ have not yet been found. These results highlight the flexibility of algorithms to incorporate additional constraints. In principle, multiple constraints can be used, the only modifications would be to add additional Gaussian process regressors and incorporate them in the calculation of the safe set. 

 \begin{figure}[t!]
    \centering
    \psfrag{-1000}[][]{\scriptsize -1000}
    \psfrag{0}[][]{\scriptsize 0~}
    \psfrag{0.5}[][]{\scriptsize 0.5~}
    \psfrag{-0.5}[][]{\scriptsize -0.5~}
    \psfrag{-1.5}[][]{\scriptsize ~~-1.5}
    \psfrag{1}[][]{\scriptsize 1~}
    \psfrag{1.5}[][]{\scriptsize 1.5~}
    \psfrag{2.5}[][]{\scriptsize 2.5~}   
    \psfrag{2}[][]{\scriptsize 2}
    \psfrag{4}[][]{\scriptsize 4~~}
    \psfrag{500}[][]{\scriptsize 500~}
    \psfrag{1000}[][]{\scriptsize 1000~}
    \psfrag{time (s)}[][]{\scriptsize time (s)}
    \psfrag{1500}[][]{\scriptsize 1500~}
    \psfrag{2000}[][]{\scriptsize 2000~}
    \psfrag{gain}[][]{\scriptsize $n$ (-)}
    \psfrag{Iteration1}[][]{\scriptsize $\iter = 1$}
    \psfrag{Iteration5}[][]{\scriptsize $\iter = 5$}
    \psfrag{Iteration10}[][]{\scriptsize $\iter = 10$}
    \psfrag{Iteration40}[][]{\scriptsize $\iter = 40$}
    \psfrag{rmsPressure}[][t]{\scriptsize O (Pa)}
    \psfrag{rmsVoltage}[][t]{\scriptsize C (V)}
    \includegraphics[trim=-0.3cm 0.05cm 0cm 0cm,clip,width=0.95\textwidth]{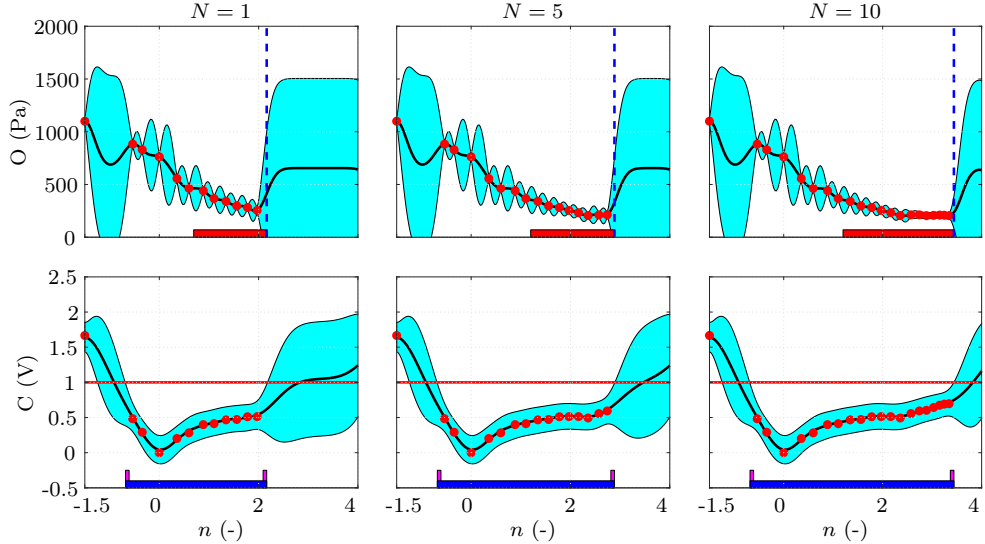}

    \caption[StageOpt Optimisation 1D]{Controller $n$ optimization with StageOpt algorithm after 1, 5, and 10 iterations. \tr{Note that $\iter$ displayed has been restarted to 1 once the additional constrain was activated after the 11th iteration.} The time delay, $\tau$, is fixed at 1.5~ms. (\bluesquare): safe set, (\redsquare): minimizer, (\magsquare): expander, (\blackline): mean prediction, (\bluedashed): next evaluation point, (\redcircle): measurement points, (\redline): safety constraint. The StageOpt algorithm was activated after ten iterations of SafeOpt. The algorithm always chooses the minimum of the lower confidence bound of the objective function for the next evaluation point.}
    \label{fig:n_stageOpt}
\end{figure}

The application of stageOpt algorithm is shown in figure \ref{fig:n_stageOpt}. Similarly to the shrinking algorithm, the first 10 iterations employ the safeOpt algorithm and then the acquisition function is switched to minimum lower confidence bound of the objective function as described in eq. \ref{eq:nextpointStageOpt}. After the switch of the acquistion function, it is clear that the algorithm always evaluates points with the lowest lower confidence bound, thereby, points with high pulsation values are not evaluated.  

\begin{figure}[t!]
    \centering
    \psfrag{250}[][]{\scriptsize 250}
    \psfrag{1}[][]{\scriptsize 1~}
    \psfrag{0.5}[][]{\scriptsize 0.5~}
    \psfrag{0.25}[][]{\scriptsize 0.25~}
    \psfrag{0.75}[][]{\scriptsize 0.75~}
    \psfrag{2}[][]{\scriptsize 2}
    \psfrag{4}[][]{\scriptsize 4}
    \psfrag{6}[][]{\scriptsize 6}
    \psfrag{250}[][]{\scriptsize 250~}
    \psfrag{750}[][]{\scriptsize 750~}
    \psfrag{500}[][]{\scriptsize 500~}
    \psfrag{1000}[][]{\scriptsize 1000~}
    \psfrag{time (s)}[][]{\scriptsize time (s)}
    \psfrag{1500}[][]{\scriptsize 1500~}
    \psfrag{2000}[][]{\scriptsize 2000~}
    \psfrag{gain}[b][]{\scriptsize $n$ (-)}
    \psfrag{Iterations}[][]{\scriptsize $\iter$}
    \psfrag{rmsPressure}[][t]{\scriptsize $O$ (Pa)}
    \psfrag{rmsVolt}[][]{\scriptsize $C$ (V)}
    \psfrag{SafeOptttt}[][]{\scriptsize SafeOpt}
    \psfrag{StageOptttt}[][]{\scriptsize StageOpt}
    \psfrag{Shrinking}[][]{\scriptsize ~Shrinking}
    \psfrag{SafetyCrit}[][]{\scriptsize Safety}
    \psfrag{rmsVolt}[][]{\scriptsize C (V)}
    \includegraphics[trim=0cm 0.5cm 0cm 0cm,clip,width=.5\textwidth]{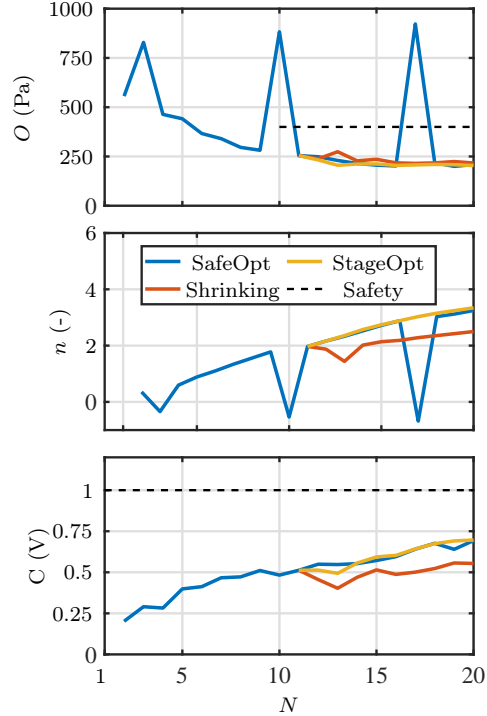}

    \caption{The evolution of the objective function, gain, and the constraint function with respect to the number of iterations. The safety constraint on the objective function is only applied to the shrinking algorithm. The first iteration contains no evaluation point.}
    \label{fig:1D_param_evol}
\end{figure}

The evolution of the values of the objective function, the gain and the constraint function, with the application of the three algorithms, are shown in figure \ref{fig:1D_param_evol}. Note that the safe constraints are indicated by the dashed line, and the constraint on the objective function is only applicable to the shrinking algorithm. As seen, in the $18^{\tr{\mathrm{th}}}$ iteration, the safeOpt algorithm evaluates points with high pulsation. On the contrary, both the stageOpt and shrinking algorithms do not have this behavior. Due to the removal of the expander, the shrinking algorithm slowly expands the safe set, as can be seen by the slowly increasing gain after the $11^{\tr{\mathrm{th}}}$ iteration. The choice of algorithms might depend on the system of interest and the preference of the user. 

 \begin{figure}[t!]
    \centering
    \psfrag{-1000}[][]{\scriptsize -1000}
    \psfrag{0}[][]{\scriptsize 0~}
    \psfrag{0.5}[][]{\scriptsize 0.5~}
    \psfrag{-0.5}[][]{\scriptsize -0.5~}
    \psfrag{-1.5}[][]{\scriptsize ~~-1.5}
    \psfrag{1}[][]{\scriptsize 1~}
    \psfrag{1.5}[][]{\scriptsize 1.5~}
    \psfrag{2.5}[][]{\scriptsize 2.5~}
    \psfrag{2}[][]{\scriptsize 2}
    \psfrag{0.5}[][]{\scriptsize 0.5}
    \psfrag{4}[][]{\scriptsize 4~~}
    \psfrag{500}[][]{\scriptsize 500~}
    \psfrag{1000}[][]{\scriptsize 1000~}
    \psfrag{time (s)}[][]{\scriptsize time (s)}
    \psfrag{1500}[][]{\scriptsize 1500~}
    \psfrag{2000}[][]{\scriptsize 2000~}
    \psfrag{gain}[][]{\scriptsize $n$ (-)}
    \psfrag{Iter1}[][]{\scriptsize $\iter = 1$}
    \psfrag{Iter5}[][]{\scriptsize $\iter = 5$}
    \psfrag{Iter10}[][]{\scriptsize $\iter = 10$}
    \psfrag{Iteration40}[][]{\scriptsize $\iter = 40$}
    \psfrag{rmsPressure}[][t]{\scriptsize O (Pa)}
    \psfrag{rmsVoltage}[][t]{\scriptsize C (V)}
    \includegraphics[trim=0cm 0.05cm 0cm 0cm,clip,width=0.95\textwidth]{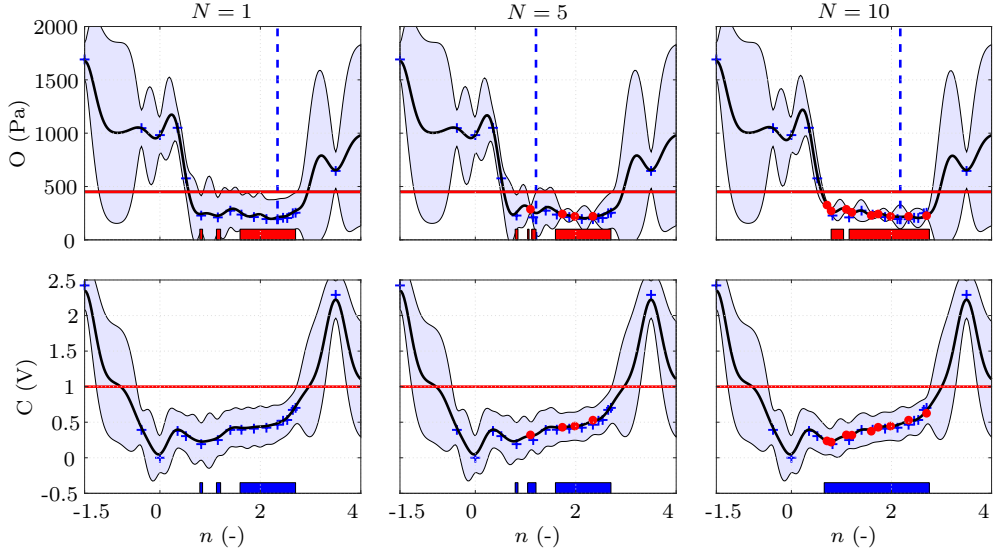}

   \caption[Context Optimisation]{Context with the shrinking algorithm after 1, 5, and 10 iterations. The algorithm makes use of the information obtained at OP2 to optimize the $n$ at OP1. The time delay, $\tau$, is fixed at 1.5~ms. (\bluesquare): safe set, (\redsquare): minimizer, (\blackline): mean prediction, (\bluedashed): next evaluation point, (\textcolor{blue}{+}): data points from OP2, (\redcircle): evaluation points, (\redline): safety constraint.}
    \label{fig:shrinkingcontext1D}

\end{figure}

The Bayesian context is applied to transfer the knowledge obtained from OP2 to OP1. The safeOpt algorithm was first performed on OP2 for 15 iterations. Afterwards, the combustor is operated to OP1 and by adding a context variable, which in this case the equivalence ratio, $\mathbf{z}=\phi$, the information from OP2 can be transferred to OP1. The length scale for the kernel $k_{\phi}$ is set to 0.1, which implies that a 0.1 difference in the equivalence ratio will lead to completely different behaviors. Figure \ref{fig:shrinkingcontext1D} shows the results of the Bayesian context. The shrinking algorithm is used, and the threshold values are the same as in Figure \ref{fig:n_shrinking}. As seen, in the first iteration, the uncertainties in the objective and constraint function are scaled up. The algorithms are now more uncertain about the information collected from OP2, which is indicated by the (\textcolor{blue}{+}) sign. This is due to the introduction of $k_\phi$ in the kernel function of both objective and constraint functions. Because of the possibility of transferring the information, the algorithm is now more sample efficient as the points with low pulsation values are now known. This is a clear advantage as small changes in operating conditions of the combustor will not require the algorithm to restart again from zero.

 \begin{figure}[t]
    \centering
    \psfrag{-1000}[][]{\scriptsize -1000}
    \psfrag{0}[][]{\scriptsize 0}
    \psfrag{1}[][]{\scriptsize 1}
    \psfrag{2}[][]{\scriptsize 2}
    \psfrag{3}[][]{\scriptsize 3}
    \psfrag{6}[][]{\scriptsize 6}
    \psfrag{-1.5}[][]{\scriptsize ~~-1.5}
    \psfrag{1}[][]{\scriptsize 1~}
    \psfrag{1.5}[][]{\scriptsize 1.5~}
    \psfrag{2.5}[][]{\scriptsize 2.5~}
    \psfrag{0.6}[][]{\scriptsize 0.6}
    \psfrag{1.2}[][]{\scriptsize 1.2}
    \psfrag{1.8}[][]{\scriptsize 1.8}
    \psfrag{200}[][]{\scriptsize 200}
    \psfrag{600}[][]{\scriptsize 600}
    \psfrag{1000}[][]{\scriptsize 1000}
    \psfrag{Gain}[][]{\scriptsize $n~(-)$}
    \psfrag{1400}[][]{\scriptsize 1400}
    \psfrag{2000}[][]{\scriptsize 2000~}
    \psfrag{gain}[][]{\scriptsize gain (-)}
    \psfrag{Pa}[][]{\scriptsize Pa}
    \psfrag{V}[][]{\scriptsize V}
    \psfrag{iter1}[][]{\scriptsize $\iter = 1$}
    \psfrag{iter30}[][]{\scriptsize $\iter = 30$}
    \psfrag{iter15}[][]{\scriptsize $\iter = 15$} 
    \psfrag{PhaseShift}[][t]{\scriptsize $\tau~(ms)$}
    \includegraphics[trim=-0.3cm 0.02cm 0cm 0cm,clip,width=0.94\textwidth] {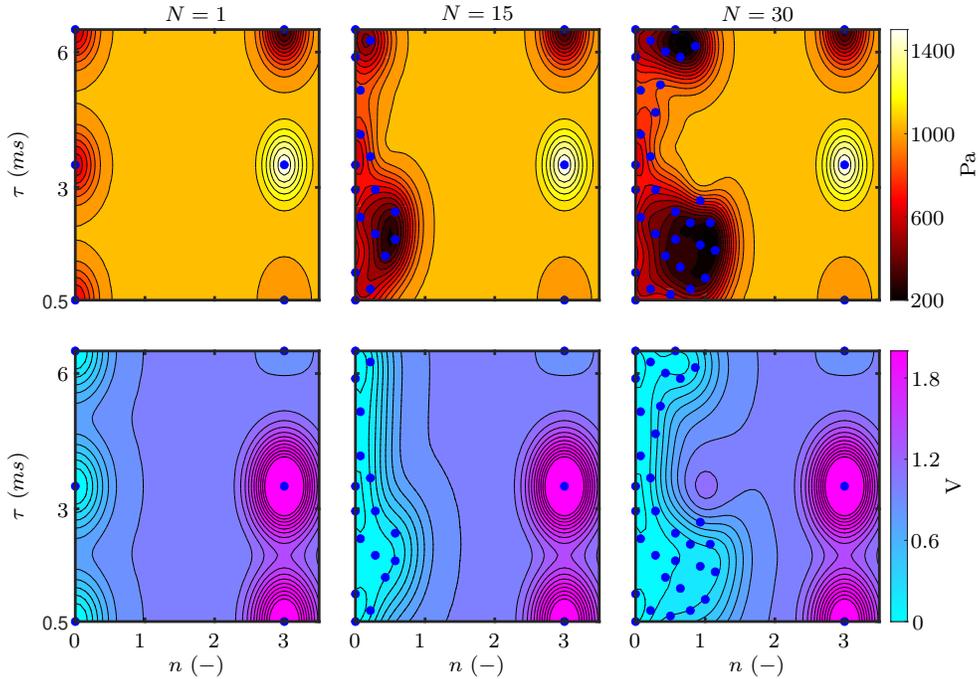}

    \caption[SafeOptimization 2D]{$n-\tau$ optimization with safeOpt algorithm after 1, 15, and 30 iterations. The mean prediction of the surface map of the objective function (top) and constraint function (bottom). (\bluecircle): evaluation points. (OP1) }
    \label{fig:n_tau_safeOpt}

\end{figure}

The framework is now extended to optimize both the $n$ and $\tau$ parameters. The control parameter space is discretized uniformly with 50$\times$50 grid points, with $n$ ranging from 0 to 3.5 and $\tau$ ranging from 0.5 to 6~ms, $\mathcal{P} \subset \mathbb{R}^{[0~3.5]\times[0.5~3.5]}$.  As a result of the increased number of parameters, a greater number of initial points is required, similar to the case of the numerical simulation. Since there is only one unstable mode observed in the experiment, only six initial points are given. Figure \ref{fig:n_tau_safeOpt} shows the mean prediction of the objective function $\mathbf{\mu_o}$ and the constraint function $\mathbf{\mu_c}$. The initial points are shown in the first column of Figure \ref{fig:n_tau_safeOpt}. The safety threshold values $T$ and $T_o$ are the same as those in the one parameter optimization case. The total number of iterations is set at 30. For the first 15 iterations, the algorithm safely explores the parameter space and finds the two regions with low pulsations, as shown in the third column of Figure \ref{fig:n_tau_safeOpt}. The first safe region is located on the bottom left region, and the second one is on the top left region. This is possible because for a single-mode instability at a frequency of $f_0$, two different time delays $\tau_2$ and $\tau_1$ that are related as $\tau_2 = \tau_1 + 1/f_{0}$ would give similar results because the phase of the controller is the same for these two delays. The algorithm evaluates more points in the bottom left region and the best evaluated point is located at $(n,\tau) = (1.25,1.7)$.

\begin{figure}[t!]
    \centering
    \psfrag{-1000}[][]{\scriptsize -1000}
    \psfrag{0}[][]{\scriptsize 0}
    \psfrag{1}[][]{\scriptsize 1}
    \psfrag{2}[][]{\scriptsize 2}
    \psfrag{3}[][]{\scriptsize 3}
    \psfrag{6}[][]{\scriptsize 6}
    \psfrag{-1.5}[][]{\scriptsize ~~-1.5}
    \psfrag{1}[][]{\scriptsize 1~}
    \psfrag{1.5}[][]{\scriptsize 1.5~}
    \psfrag{2.5}[][]{\scriptsize 2.5~}
    \psfrag{0.6}[][]{\scriptsize 0.6}
    \psfrag{1.2}[][]{\scriptsize 1.2}
    \psfrag{1.8}[][]{\scriptsize 1.8}
    \psfrag{200}[][]{\scriptsize 200}
    \psfrag{600}[][]{\scriptsize 600}
    \psfrag{1000}[][]{\scriptsize 1000}
    \psfrag{Gain}[][]{\scriptsize $n~(-)$}
    \psfrag{1400}[][]{\scriptsize 1400}
    \psfrag{2000}[][]{\scriptsize 2000~}
    \psfrag{gain}[][]{\scriptsize gain (-)}
    \psfrag{Pa}[][]{\scriptsize Pa}
    \psfrag{V}[][]{\scriptsize V}
    \psfrag{iter1}[][]{\scriptsize $\iter = 1$}
    \psfrag{iter15}[][]{\scriptsize $\iter = 15$}
    \psfrag{iter7}[][]{\scriptsize $\iter = 7$}
    \psfrag{PhaseShift}[][t]{\scriptsize $\tau~(ms)$}
    
    \includegraphics[trim=-0.3cm 0.02cm 0cm 0cm,clip,width=0.94\textwidth]{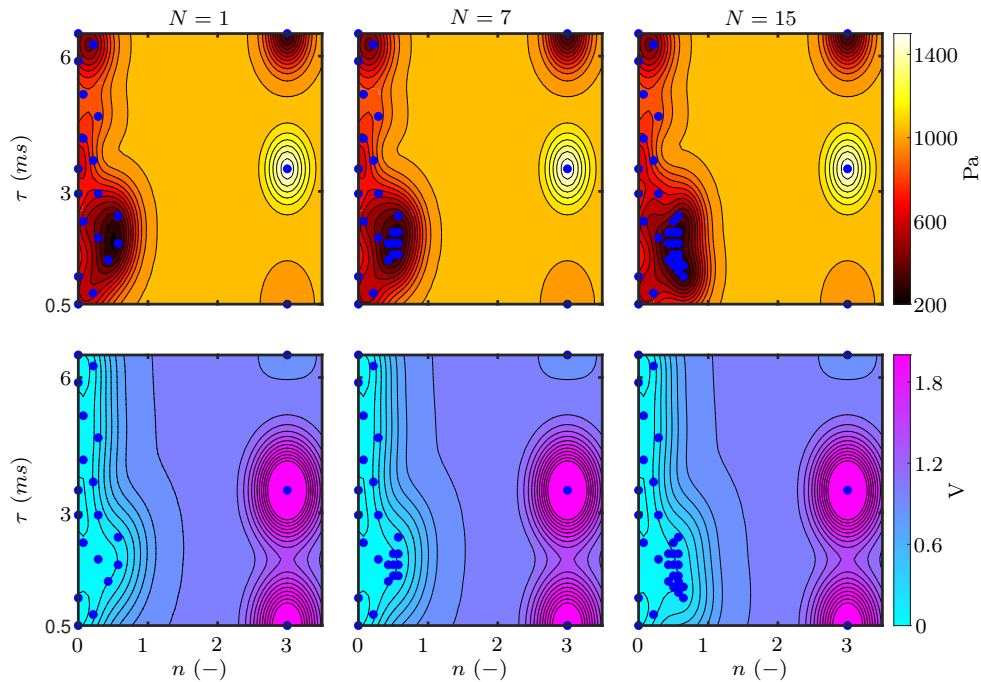}

    \caption[Shrinking 2D]{$n-\tau$ optimization with Shrinking algorithm after 1, 7, and 15 iterations. The first iteration is taken from the $\mathrm{15^{th}}$ iteration of the safeOpt. The mean prediction of the surface map of the objective function (top) and constraint function (bottom). (\bluecircle): evaluation points. The expansion of the parameter space becomes more restricted due to an additional constraint on the objective function.}
    \label{fig:n_tau_shrinking}

\end{figure}

The results obtained through the use of the shrinking algorithm to optimize the two parameters are presented in figure \ref{fig:n_tau_shrinking}. In the initial 15 iterations, the algorithm employs the safeOpt approach, followed by the introduction of an additional constraint applied to the objective function. Similarly to the previous scenario that involves single-parameter optimization, the expander is eliminated. This change results in a more constrained focus on evaluating points in the bottom-left region. Due to the absence of expanders, the assessed points in the lower left region are positioned closely together, significantly limiting the algorithm's ability to expand. Nevertheless, the safe-set region continues to grow over the course of the iterations. It is important to note that the majority of the evaluated points exhibit thermoacoustic stability. However, they tend to have slightly higher root mean square (rms) values compared to the best-evaluated point obtained with the safeOpt algorithm. One potential approach to enable more aggressive exploration would involve reintroducing the expander into the decision-making process for the next point evaluation. Nevertheless, the primary objective here is to illustrate that there can be drawbacks when the exploratory aspect is curtailed.

\begin{figure}[t!]
    \centering
        \psfrag{-1000}[][]{\scriptsize -1000}
    \psfrag{0}[][]{\scriptsize 0}
    \psfrag{1}[][]{\scriptsize 1}
    \psfrag{2}[][]{\scriptsize 2}
    \psfrag{3}[][]{\scriptsize 3}
    \psfrag{6}[][]{\scriptsize 6}
    \psfrag{-1.5}[][]{\scriptsize ~~-1.5}
    \psfrag{1}[][]{\scriptsize 1~}
    \psfrag{1.5}[][]{\scriptsize 1.5~}
    \psfrag{2.5}[][]{\scriptsize 2.5~}
    \psfrag{0.6}[][]{\scriptsize 0.6}
    \psfrag{1.2}[][]{\scriptsize 1.2}
    \psfrag{1.8}[][]{\scriptsize 1.8}
    \psfrag{200}[][]{\scriptsize 200}
    \psfrag{600}[][]{\scriptsize 600}
    \psfrag{1000}[][]{\scriptsize 1000}
    \psfrag{Gain}[][]{\scriptsize $n~(-)$}
    \psfrag{1400}[][]{\scriptsize 1400}
    \psfrag{2000}[][]{\scriptsize 2000~}
    \psfrag{gain}[][]{\scriptsize gain (-)}
    \psfrag{Pa}[][]{\scriptsize Pa}
    \psfrag{V}[][]{\scriptsize V}
    \psfrag{iter1}[][]{\scriptsize $\iter = 1$}
    \psfrag{iter7}[][]{\scriptsize $\iter = 7$}
    \psfrag{iter15}[][]{\scriptsize $\iter = 15$}
    \psfrag{PhaseShift}[][t]{\scriptsize $\tau~(ms)$}
    
    \includegraphics[trim=-0.3cm 0.02cm 0cm 0cm,clip,width=0.94\textwidth]{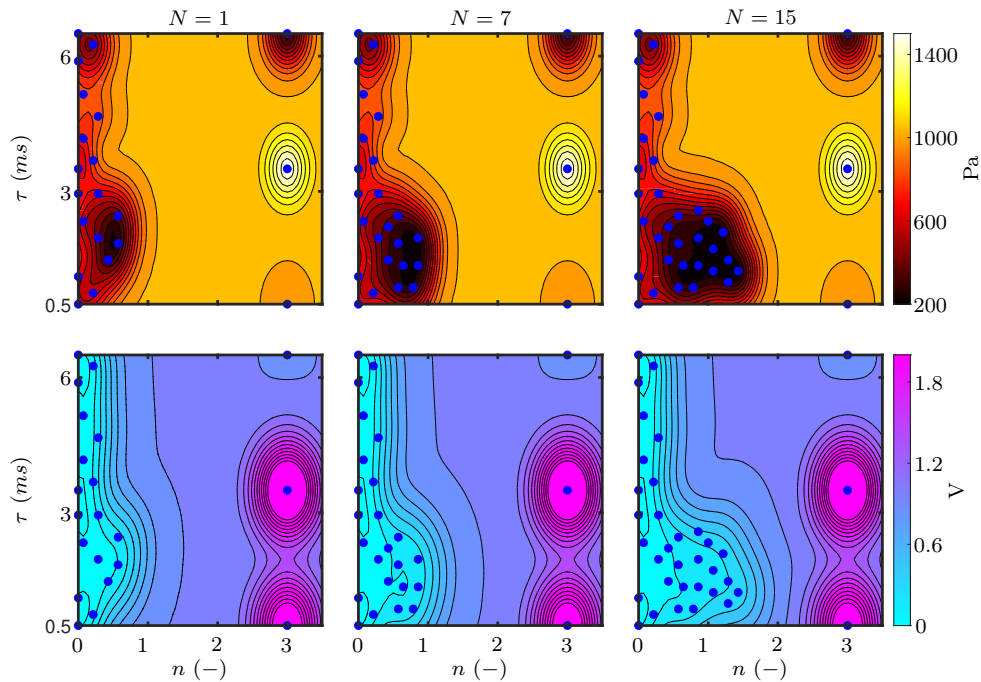}

    \caption[StageOpt 2D]{$n-\tau$ optimization with stageOpt algorithm after 1, 7, and 15 iterations. The first iteration is taken from the $\mathrm{15^{th}}$ iteration of the safeOpt. The mean prediction of the surface map of the objective function (top) and constraint function (bottom). (\bluecircle) evaluation points. The algorithm picks more points  on the lower left plane of the parameter space. (OP1)}
    \label{fig:n_tau_stageOpt}

\end{figure}

The results obtained using the stageOpt algorithm are depicted in figure \ref{fig:n_tau_stageOpt}. Similarly to the shrinking algorithm, the initial 15 iterations employ the safeOpt algorithm for a cautious exploration of the parameter space. Similarly to the scenario involving single-parameter optimization, following the transition to a different acquisition function, the algorithm refrains from assessing points associated with high pulsation values. The stageOpt algorithm, however, evaluates regions with $\tau$ values akin to those targeted by the safeOpt algorithm, falling within the range of 0.5 to 2 milliseconds. Notably, the stageOpt algorithm extends its evaluations to higher gain values. This expanded exploration of the parameter space is facilitated because the alternative safe region, characterized by low pulsation values at $\tau$ around 6 milliseconds, remains unexplored. As a result, with an equal number of total iterations, the stageOpt algorithm efficiently allocates more iterations to the assessment of points in the lower left region of the domain.

\begin{figure}[t!]
    \centering
    \psfrag{250}[][]{\scriptsize 250}
    \psfrag{1}[][]{\scriptsize 1~}
    \psfrag{0}[][]{\scriptsize 0}
    \psfrag{15}[][]{\scriptsize 15}
    \psfrag{30}[][]{\scriptsize 30}
    \psfrag{150}[][]{\scriptsize 0}
    \psfrag{675}[][]{\scriptsize 675}
    \psfrag{1200}[][]{\scriptsize 1200}
    \psfrag{3.5}[][]{\scriptsize 3.5}
    \psfrag{7}[][]{\scriptsize 7}
    \psfrag{0.5}[][]{\scriptsize 0.5}
    \psfrag{1.5}[][]{\scriptsize 1.5}
    \psfrag{2}[][]{\scriptsize 2}
    \psfrag{4}[][]{\scriptsize 4}
    \psfrag{6}[][]{\scriptsize 6}
    \psfrag{250}[][]{\scriptsize 250~}
    \psfrag{750}[][]{\scriptsize 750~}
    \psfrag{500}[][]{\scriptsize 500~}
    \psfrag{1000}[][]{\scriptsize 1000~}
    \psfrag{time (s)}[][]{\scriptsize time (s)}
    \psfrag{1500}[][]{\scriptsize 1500~}
    \psfrag{2000}[][]{\scriptsize 2000~}
    \psfrag{gain}[b][]{\scriptsize $n$ (-)}
    \psfrag{taums}[][]{\scriptsize $\tau$ (ms)}
    \psfrag{Constraint}[b][]{\scriptsize $C$ (V)}
    \psfrag{safeOpt1234}[][]{\scriptsize safeOpt}
    \psfrag{stageOpt1234}[][]{\scriptsize stageOpt}
    \psfrag{shrinking1234}[][]{\scriptsize shrinking}
    \psfrag{safety1234}[][]{\scriptsize safety}
    \psfrag{iterations}[][]{\scriptsize $\iter$}
    \psfrag{Objective}[][t]{\scriptsize $O$ (Pa)}
    \psfrag{rmsVolt}[][]{\scriptsize $C$ (V)}
    \psfrag{SafeOptttt}[][]{\scriptsize SafeOpt}
    \psfrag{StageOptttt}[][]{\scriptsize StageOpt}
    \psfrag{Shrinking}[][]{\scriptsize ~Shrinking}
    \psfrag{SafetyCrit}[][]{\scriptsize Safety}
    \psfrag{rmsVolt}[][]{\scriptsize C (V)}
    \includegraphics[trim=-0.3cm 0cm 0.0cm 0cm,clip,width=1\textwidth]{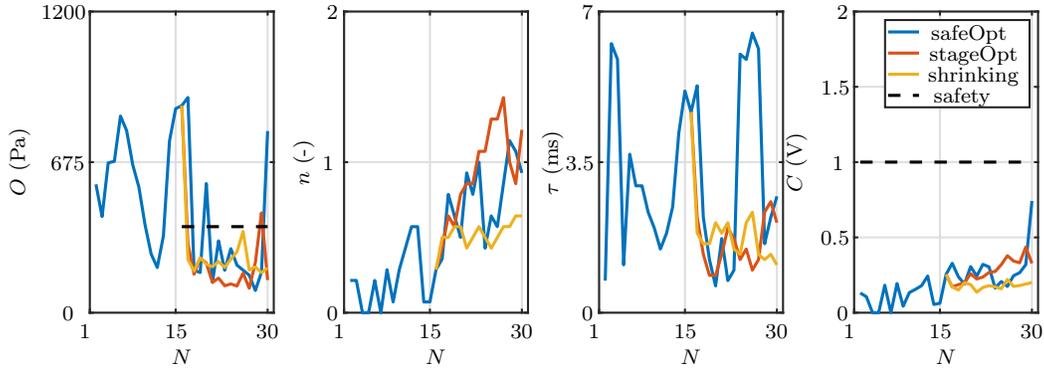}

    \caption{The evolution of the objective function, gain, and the constraint function with respect to the number of iterations. The safety constraint on the objective function is only applied to the shrinking algorithm. The first iteration does not contain an evaluation point.}
    \label{fig:2D_param_evol}
\end{figure}

The figure \ref{fig:2D_param_evol} presents the progression of various parameters, including the objective function, gain $n$, delay $\tau$, and the constraint function. It's evident from the graph that the safeOpt algorithm explores a region where $\tau$ is approximately 6.2 milliseconds, from iteration $\iter$ 22 to 27. In contrast, both the stageOpt and shrinking algorithms consistently assess regions where $\tau$ is less than 2.5 milliseconds. This difference arises because, as mentioned earlier, within the first 15 iterations of the safeOpt algorithm, the region characterized by low objective function values with $\tau$ around 6.2 milliseconds has not yet been discovered. Consequently, both the stageOpt and shrinking algorithms remain unaware of this particular region. However, it's worth noting that the region with $\tau$ around 6.2 milliseconds is expected to yield performance similar to that of $\tau$ around 1.5 milliseconds. As a result, the stageOpt and shrinking algorithms do not suffer any disadvantages in terms of failing to identify the global minimum, given that this region eventually provides equivalent performance.

In different scenarios, increasing the value of $N_{s}$ in both Algorithm \ref{alg:shrinking} and Algorithm \ref{alg:stageOpt} could provide a more comprehensive overview of the parameter space for the algorithm. However, it is important to emphasize that the choice of $N_{s}$ is problem-specific and contingent on the particular case at hand. Consequently, users should establish their expectations before initiating the optimization process.

\begin{figure}[t!]
    \centering
    \psfrag{-1000}[][]{\scriptsize -1000}
    \psfrag{0}[][]{\scriptsize 0~}
    \psfrag{1}[][]{\scriptsize 1}
    \psfrag{2}[][]{\scriptsize 2}
    \psfrag{3}[][]{\scriptsize 3~}
    \psfrag{6}[][]{\scriptsize 6~}
    \psfrag{-1.5}[][]{\scriptsize ~~-1.5}
    \psfrag{1}[][]{\scriptsize 1~}
    \psfrag{1.5}[][]{\scriptsize 1.5~}
    \psfrag{2.5}[][]{\scriptsize 2.5~}
    \psfrag{0.6}[][]{\scriptsize 0.6}
    \psfrag{1.2}[][]{\scriptsize 1.2}
    \psfrag{1.8}[][]{\scriptsize 1.8}
    \psfrag{200}[][]{\scriptsize ~~200}
    \psfrag{800}[][]{\scriptsize ~~800}
    \psfrag{1000}[][]{\scriptsize 1000}
    \psfrag{Gain}[][]{\scriptsize $n~(-)$}
    \psfrag{1400}[][]{\scriptsize ~~1400}
    \psfrag{2000}[][]{\scriptsize 2000~}
    \psfrag{gain}[][]{\scriptsize $n$ (-)}
    \psfrag{Pa}[][]{\scriptsize Pa}
    \psfrag{V}[][]{\scriptsize V}
    \psfrag{iter1OC1}[][]{\scriptsize $\iter = 1$, OP2}
    \psfrag{iter15OC1}[][]{\scriptsize $\iter = 15$, OP2}
    \psfrag{iter30OC1}[][]{\scriptsize $\iter = 30$, OP2}
    \psfrag{iter1OC2}[][]{\scriptsize $\iter = 1$, OP1}
    \psfrag{iter4OC2}[][]{\scriptsize $\iter = 4$, OP1}
    \psfrag{iter8OC2}[][]{\scriptsize $\iter = 8$, OP1}
    \psfrag{iter1OC3}[][]{\scriptsize $\iter = 1$, OP3}
    \psfrag{iter4OC3}[][]{\scriptsize $\iter = 4$, OP3}
    \psfrag{iter8OC3}[][]{\scriptsize $\iter = 8$, OP3}
    \psfrag{iter7}[][]{\scriptsize $\iter = 7$}
    \psfrag{phaseshift}[b][t]{\scriptsize $\tau$~(ms)}
    \includegraphics[trim=2.1cm 1.1cm 0.7cm 1.cm,clip,width=0.94\textwidth]{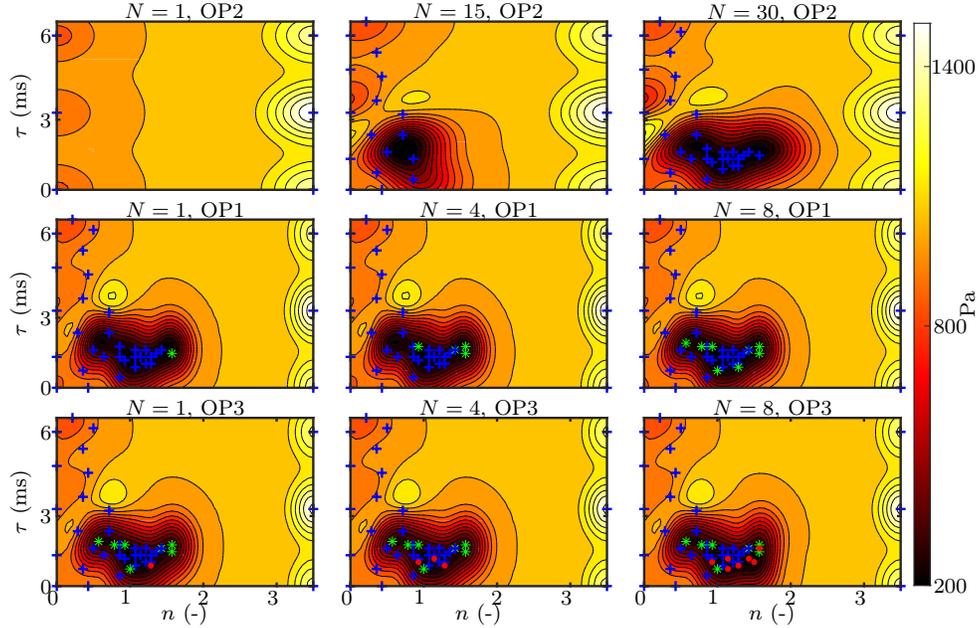}

    \caption[2D context]{$n-\tau$ optimization with context. The algorithm uses the context variable to transfer the knowledge between each operating point. The mean predictions of the objective function from OP2, OP1 and OP3 are displayed in the second, first, and third rows, respectively. \textcolor{blue}{+}: evaluation points from OP1, \textcolor{green}{*}: evaluation points from OP2,(\redcircle): evaluation points from OP3.}
    \label{fig:n_tau_shrinking_context}

\end{figure}

\begin{figure}[t!]
    \centering
    \psfrag{-1000}[][]{\scriptsize -1000}
    \psfrag{0}[][]{\scriptsize 0}
    \psfrag{1}[][]{\scriptsize 1}
    \psfrag{2}[][]{\scriptsize 2}
    \psfrag{3}[][]{\scriptsize 3}
    \psfrag{6}[][]{\scriptsize 6}
    \psfrag{-1.5}[][]{\scriptsize ~~-1.5}
    \psfrag{1}[][]{\scriptsize 1~}
    \psfrag{1.5}[][]{\scriptsize 1.5~}
    \psfrag{2.5}[][]{\scriptsize 2.5~}
    \psfrag{0.5}[][]{\scriptsize 0.5}
    \psfrag{1.2}[][]{\scriptsize 1.2}
    \psfrag{1.8}[][]{\scriptsize 1.8}
    \psfrag{200}[][]{\scriptsize ~200}
    \psfrag{800}[][]{\scriptsize ~800}
    \psfrag{1000}[][]{\scriptsize 1000}
    \psfrag{Gain}[][]{\scriptsize $n~(-)$}
    \psfrag{1400}[][]{\scriptsize 1400}
    \psfrag{2000}[][]{\scriptsize 2000~}
    \psfrag{gain}[][]{\scriptsize gain (-)}
    \psfrag{Pa}[l][c][1][270]{\scriptsize $\frac{\sigma}{\sigma_{max}}$}
    \psfrag{V}[][]{\scriptsize V}
    \psfrag{iter1OC1}[][]{\scriptsize $\iter = 1$, OP2}
    \psfrag{iter15OC1}[][]{\scriptsize $\iter = 15$, OP2}
    \psfrag{iter30OC1}[][]{\scriptsize $\iter = 30$, OP2}
    \psfrag{iter1OC2}[][]{\scriptsize $\iter = 1$, OP1}
    \psfrag{iter4OC2}[][]{\scriptsize $\iter = 4$, OP1}
    \psfrag{iter8OC2}[][]{\scriptsize $\iter = 8$, OP1}
    \psfrag{iter1OC3}[][]{\scriptsize $\iter = 1$, OP3}
    \psfrag{iter4OC3}[][]{\scriptsize $\iter = 4$, OP3}
    \psfrag{iter8OC3}[][]{\scriptsize $\iter = 8$, OP3}
    \psfrag{iter7}[][]{\scriptsize $\iter = 7$}
    \psfrag{PhaseShift}[][t]{\scriptsize $\tau~(ms)$}
    \includegraphics[trim=-0.3cm 0.02cm 0cm 0cm,clip,width=0.94\textwidth]{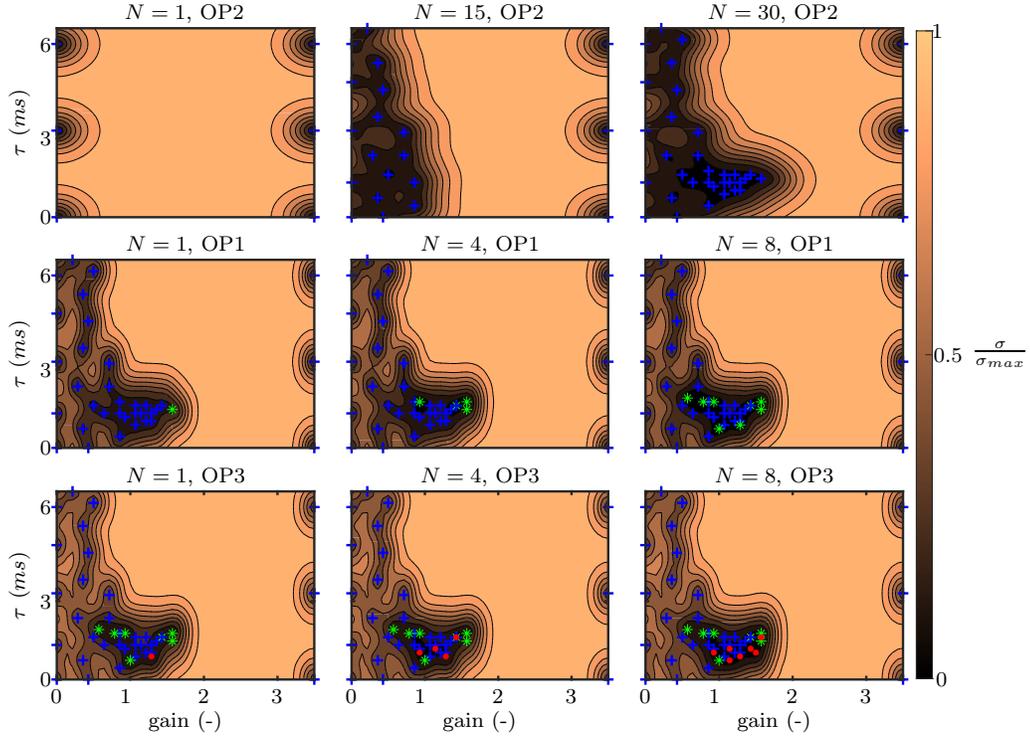}

    \caption[Uncertainty Map]{Map of the normalized uncertainty of the objective function for $n-\tau$ optimization with context. The algorithm uses the context variable to transfer the knowledge between each operating point. The mean predictions of the objective function from OP1, OP2, and OP3 are displayed in the first, second, and third rows, respectively. \textcolor{blue}{+}: evaluation points from OP1, \textcolor{green}{*}: evaluation points from OP2,(\redcircle): evaluation points from OP3.}
    \label{fig:n_tau_shrinking_context_uncertainty}

\end{figure}

Bayesian context is applied for two-parameter optimization, and three operating conditions are considered. The shrinking algorithm is applied first to OP2 with the number of iterations equal to 30 and $N_s$ = 15, the expander is used to acquire the next point evaluation. Subsequently, the information is transferred to OP1, then the shrinking algorithm is applied with the number of iterations equals to 8 with $N_s$ = 1. Subsequently, the information obtained from both OP2 and OP1 is transferred to OP3, and the shrinking algorithm with the total number of iterations equals to 8 and $N_s$ = 1 is applied. The kernel parameters of the objective and constraint functions, as well as the safety threshold values, are the same as in the previous case. The mean prediction of the objective function for the three operating conditions is shown in figure \ref{fig:n_tau_shrinking_context}, whereas the normalized uncertainty of the objective function is shown in figure \ref{fig:n_tau_shrinking_context_uncertainty}. 

As can be seen from the first row of Figure \ref{fig:n_tau_shrinking_context}, the algorithm can safely explore the parameter space and the evaluation points are less closely packed compared to the result in figure \ref{fig:n_tau_shrinking}. This is due to the inclusion of the expander in the shrinking algorithm which allows the algorithm to explore more aggressively while still respecting the safety threshold in both the objective and constraint functions. Because the frequency of the instability at OP2 is lower than at OP1, the other region with low pulsation values as in figure \ref{fig:n_tau_safeOpt} lies outside the domain. After completing 30 iterations at OP2, the information is carried out to OP1 with the Bayesian context. As clearly seen from the second row of figure \ref{fig:n_tau_shrinking_context}, the algorithm evaluates points in the vicinity where the previous operating point exhibits low pulsation and is safe.

Following the completion of eight iterations at OP2, the combustor transitions to OP3, with information from both OP2 and OP1 being carried over to OP3. In the graph displayed in figure \ref{fig:n_tau_shrinking_context_uncertainty}, it becomes evident that, at OP3, the algorithm exhibits increased uncertainty regarding the information derived from OP2 in comparison to that from OP1. This heightened uncertainty is attributed to the fact that the context variable of OP1 is closer to OP2 than it is to OP3. A similar pattern emerges at OP3, where the algorithm now recognizes the safe region with low pulsations and consistently evaluates points within this low pulsation region. While it's possible that eight iterations may be insufficient for the algorithm to thoroughly explore the parameter space, the primary aim here is to illustrate that with the aid of Bayesian context, direct access to a safe region with low pulsation values can be achieved without the necessity of restarting the algorithm from scratch.

\subsubsection{Sequential combustor with NRPD}

In the sequential combustor configuration, the sequential flame is stabilized with the help of NRPD. The objective function is the NO emission in parts per million by volume, dry (ppmvd), and the constraint function is the RMS pressure in the sequential combustion chamber, measured by the mic. 2 in figure \ref{fig:ExpSetupSequential}. The two control parameters are the voltage output of the high voltage generator and the allocation of power between both flames. It is important to note that the overall thermal power of the flames remains at a constant value of 73.4 kW, and the mass flow rates for both the first-stage air and dilution air are also maintained at a consistent level. The term "power splitting" is defined as the percentage representing the proportion of the total power directed towards the second stage flame. For example, if the power splitting is set to 60 $\%$, this means that the second stage flame thermal power is 60$\%$ of the whole thermal power, which equates to 44.04~kW. Regarding the fuel composition, it is important to note that the first stage exclusively uses 100$\%$ natural gas, while the second stage employs a mixture consisting of 10.45$\%$ hydrogen and 89.55$\%$ natural gas by mass. The emission analyzer requires a continuous operation of about 100 seconds to converge, whereas, the rms pressure is computed after recording the pressure signal for five seconds. 

\begin{figure}
    \centering
        \psfrag{aaa}[][]{\scriptsize a)~~}
    \psfrag{bbb}[][]{\scriptsize b)~~}
    \psfrag{ccc}[][]{\scriptsize c)~~}
    \psfrag{ddd}[][]{\scriptsize d)~~}
    \psfrag{eee}[][]{\scriptsize e)~~}
    \psfrag{fff}[][]{\scriptsize f)~~} 
    \psfrag{a}[][]{\scriptsize a}
    \psfrag{ b}[][]{\scriptsize b}
    \psfrag{ c}[][]{\scriptsize c}
    \psfrag{ d}[][]{\scriptsize d}
    \psfrag{ e}[][]{\scriptsize e}
    \psfrag{ f}[][]{\scriptsize f}  
    \psfrag{iter1}[][]{\scriptsize $\iter$ = 1}  
    \psfrag{iter6}[][]{\scriptsize $\iter$ = 6}  
    \psfrag{iter12}[][]{\scriptsize $\iter$ = 12}  
    \psfrag{g}[][]{\scriptsize g)}
    \psfrag{h}[][]{\scriptsize h)}
    \psfrag{PowerRatio}[][]{\scriptsize Power splitting}
     \psfrag{Voltage}[b][t]{\scriptsize $V$ (kV)}
    \psfrag{Pressure}[][]{\scriptsize P (mbar)}
    \psfrag{mixingchanqqqqqq}[][]{\scriptsize Mixing Channel }
    \psfrag{seqcomb}[][]{\scriptsize ~~~~~~~~~~~~~~Seq. Combustor }
    \psfrag{Firststageqqqqqqqqq}[][]{\scriptsize~ First combustor}
    \psfrag{ secondstage}[][]{\scriptsize ~~~~~~~~~~~~~Seq. combustor}
    \psfrag{timems}[][]{\scriptsize{$t$ (ms)}}
    \psfrag{0.4}[][]{\scriptsize{40}}
    \psfrag{0.48}[][]{\scriptsize{48}}
    \psfrag{0.56}[][]{\scriptsize{56}}
    \psfrag{0}[][]{\scriptsize{~0}}
    \psfrag{1000}[][]{\scriptsize{~1000}}
    \psfrag{1500}[][]{\scriptsize~{1500}}
    \psfrag{480}[][]{\scriptsize{480}}
    \psfrag{490}[][]{\scriptsize{490}}
    \psfrag{500}[][]{\scriptsize{~500}}
    \psfrag{7.5}[][]{\scriptsize{7.5}}
    \psfrag{8}[][]{\scriptsize{8}}
    \psfrag{8.5}[][]{\scriptsize{8.5}}
    \psfrag{9}[][]{\scriptsize{9}}
    \psfrag{5}[][]{\scriptsize{~5}}
    \psfrag{10}[][]{\scriptsize{~10}}
    \psfrag{15}[][]{\scriptsize{~15}}
    \psfrag{20}[][]{\scriptsize{~20}}
    \psfrag{25}[][]{\scriptsize{~25}}
    \psfrag{ppmvd}[t][]{\scriptsize{ppmvd}}
    \psfrag{Pa}[t][]{\scriptsize{Pa}}

    \includegraphics[trim=-0.7cm 0cm -0.4cm 0cm,clip,width=1\textwidth]{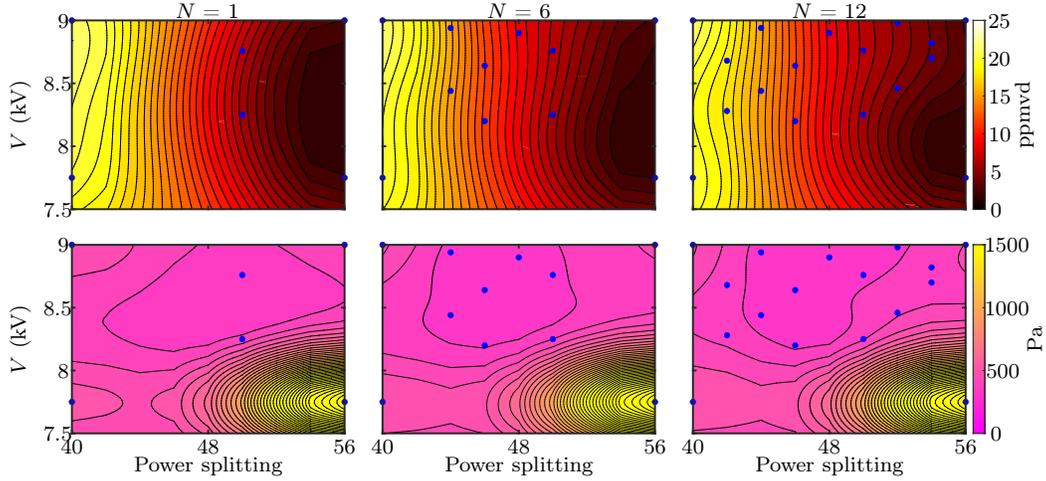}

    \caption{Bayesian optimization using the safeOpt algorithm. The objective function is the NO emissions, and the constraint function is the RMS pressure. (top) The mean prediction of the objective function. (bottom) The mean prediction of the constraint function. The blue dots show the evaluation points.}
    \label{fig:safeOpt2DNRPD}
\end{figure}

\begin{table}[t!]
\centering
\begin{tabular}{@{}ccc@{}}
\toprule
\multicolumn{1}{l}{\textbf{Hyperparameters}} & \textbf{$O$} & \textbf{$C$} \\ \midrule
$\theta$ {[}ppmvd - Pa{]}                    & 15           & 300          \\
$l_{V}$ {[}kV{]}                             & 0.5          & 0.5          \\
$l_{PS}$ {[}\%{]}                            & 2            & 2            \\
$\sigma$ {[}ppmvd - Pa{]}                    & 1.2          & 20           \\ \bottomrule
\end{tabular}
\caption{Hyperparameters of the Gaussian Process Regressors.}\label{tab:hyperparamsNRPD}
\end{table}

The control parameter space is discretized with the power splitting ranges from 40$\%$ to 56$\%$ with a step of 2$\%$, whereas the generator voltage is from 7.5~kV to 9~kV, with 100 points in between. The length scale matrix $\mathbf{L}$ is now written as: 

\begin{equation}
\mathbf{L} = \begin{bmatrix}
l_V & 0\\
0 & l_{PS},
\end{bmatrix}
\end{equation}
where $l_V$ is the length scale of the plasma generator voltage and $l_{PS}$ is the power splitting. The hyperparameters are listed in Table \ref{tab:hyperparamsNRPD}.

The safeOpt algorithm is utilized and six initial points are used to get the initial safe set for the algorithm. The initial points can be seen in the first row of figure \ref{fig:safeOpt2DNRPD}. The safety threshold for the rms pressure is set at 500 Pa. As can be seen, a power splitting of 56$\%$ with a voltage of 7.75~kV exhibits low NO emissions but high \tr{amplitude acoustic} pressure. Therefore, this point is classified as unsafe. As the iteration progresses, the algorithm safely explores the control parameter space by evaluating the region in the middle of the domain. Afterwards, the algorithm evaluates the region on the top right part of the domain. After 12 iterations, the optimization process is terminated, and the safe minimum point is found at a power split of 56$\%$ and voltage of 9~kV. As mentioned in \cite{Dharmaputra2023}, the application of NRPD \tr{can significantly} alter the \tr{mean} sequential flame position and\tr{ its heat release rate response to acoustic perturbations}, which subsequently affect the \tr{thermoacoustic} stability of the whole system. However, the main objective of this study is to optimize the control parameters solely, and the discussion about the stabilization mechanism will be investigated in a separate study. 

These results highlight the flexibility of the algorithm in optimizing different control parameters for thermoacoustic control. The control structure does not have to be affine feedback control based; as demonstrated here, the NRPD is operated in continuous fashion without feeding the pressure signal to the controller as in the previous section. Moreover, the versatility of the algorithm extends to the straightforward inclusion of additional constraint variables, such as exhaust temperature, carbon monoxide (CO) emissions, and so on. These variables can be readily incorporated as constraint functions by configuring Gaussian Process Regression with appropriately tuned kernel parameters. Furthermore, NRPD equipped with the proposed algorithms is shown to be an effective actuator for controlling both the pulsation of a sequential combustor and accessing operating conditions with low NO emissions. 

\section{Conclusions and outlook}

This study has effectively showcased the practicality of employing safe Bayesian optimization algorithms for thermoacoustic control. We have introduced and put into practice three distinct algorithms tailored specifically for thermoacoustic systems, revealing their efficacy through numerical simulations and experimental implementations. In all setups, the algorithm does not employ any model and relies on the obtained measurements to update the regressors. 

In the first phase, encompassing both numerical simulations and the first experimental setup, the adaptive optimization of feedback control parameters within a single-stage combustor utilizing loudspeaker actuation is illustrated. These proposed algorithms facilitate the transfer of knowledge between varying operating conditions, resulting in a significant acceleration of the optimization process.

In the second experimental scenario, we applied the same algorithm to a sequential combustor employing nanosecond repetitive pulsed discharges (NRPD). This setup differs significantly from the previous one in terms of combustor architecture, control actuation, and safety criteria, yet the proposed algorithms proved to be versatile and effective across these diverse contexts. 

\section*{Acknowledgments}
This project has received funding from the European Research Council (ERC) under the European Union’s Horizon 2020 research and innovation program (grant agreement No [820091]).

\section*{Author contributions}

B.D. and N.N. conceived the research project. B.D. and P.R. developed the Bayesian optimizers. B.D., P.R. and S.S. performed the experiments. All authors discussed the results. N.N. supervised the project. B.D. wrote the paper. The final version of the manuscript has been edited and approved by all the authors.

\clearpage
\bibliographystyle{elsarticle-num.bst}

\bibliography{references}
\end{document}